\documentclass[12pt,reqno,tbtags]{amsart}
\usepackage{amssymb}

\nopagebreak[3]

\numberwithin{equation}{section}

\let\Sum=\sum
\DeclareMathOperator*{\textsum}{{\textstyle\Sum}}
\def\sum{\textsum}
\let\Prod=\prod
\DeclareMathOperator*{\textprod}{{\textstyle\Prod}}
\def\prod{\textprod}
\let\Bigcap=\bigcap
\DeclareMathOperator*{\textbigcap}{{\textstyle\Bigcap}}
\def\bigcap{\textbigcap}

   \newtheorem{theorem}{Theorem}[section]
   \newtheorem{proposition}{Proposition}[section]
   \newtheorem{lemma}{Lemma}[section]
   \newtheorem{corollary}{Corollary}[section]

\theoremstyle{remark}
   \newtheorem*{remark}{\rm \textit{Remark}}
\theoremstyle{definition}

   \newtheorem{assumption}{Assumption}

\makeatletter
\newcommand{\thmskip}{\vspace{.5\baselineskip\@plus.2\baselineskip
                                    \@minus.2\baselineskip}\par\noindent}
\makeatother

\newcommand{\ebox}{~$\quad\square$}
\newcommand{\hdot}{{\dot H}}
\newcommand{\calh}{{\mathcal H}}
\newcommand{\A}{\textbf{\textit{A}}}
\newcommand{\B}{\textbf{\textit{B}}}
\newcommand{\J}{\textbf{\textit{J}}}

\newcommand{\R}{\textbf{\textit{R}}}

\newcommand{\C}{\textbf{\textit{C}}}

\newcommand{\tild}{{\tilde d}}
\newcommand{\tilq}{{\tilde q}}
\newcommand{\tilr}{{\tilde r}}
\newcommand{\tilD}{{\tilde D}}
\newcommand{\tilE}{{\tilde E}}
\newcommand{\LG}{\mathrm{L}}
\newcommand{\CG}{\mathrm{C}}

\newcommand{\embedding}{\hookrightarrow}
\newcommand{\Proof}{\par\noindent {\it Proof}. }
\newcommand{\Proofof}[1]{\par\noindent {\it Proof of #1}.}

\newcommand{\I}{1}

\def\Div{\mathop{\mathrm{div}}\nolimits}
\def\Re{\mathop{\mathrm{Re}}\nolimits}
\def\Im{\mathop{\mathrm{Im}}\nolimits}

\title[Local well-posedness for the Maxwell-Schr\"odinger equation]%
{Local well-posedness for \\ the Maxwell-Schr\"odinger equation}

\author[M. Nakamura and T. Wada]{Makoto NAKAMURA \ \ and \ \ Takeshi WADA}
\begin{document}

\maketitle

\begin{abstract}
Time local well-posedness for the Maxwell-Schr\"odinger equation in the Coulomb gauge
is studied in Sobolev spaces by the contraction mapping principle. 
The Lorentz gauge and the temporal gauge cases are also treated by the gauge transform.

\vspace{9pt}

\noindent {\bf Mathematics Subject Classification (2000):} 
%35A05, %General existence and uniqueness theorems
%35Q60, %Equations of electromagnetic theory and optics
%35L70, %Nonlinear second-order PDE of hyperbolic type
%78A25, %Electromagnetic theory, general
%81T13. %Yang-Mills and other gauge theories
%
35Q40, % Equations from quantum mechanics
35Q55, %NLS-like (nonlinear Schr\"odinger) equations
%35Q60, %Equations of electromagnetic theory and optics
35L70. 

\end{abstract}

\section{Introduction}
We consider the Maxwell-Schr\"odinger equation (\textbf{MS}):
\begin{gather}
i\partial_tu=(\calh(\A)+\phi)u,   \label{MS1} \\
-\Delta \phi-\partial_t \Div \A =\rho(u),  \label{MS2} \\
(\partial_t^2-\Delta)\A+\nabla(\partial_t\phi+\Div \A) =\J(u,\A), \label{MS3}
\end{gather}
where $(u,\phi,\A): \R^{1+3}\rightarrow \C\times\R\times\R^3$, 
$\calh(\A)=-(\nabla-i\A)^2$, $\rho(u)=|u|^2$, $\J(u,\A)=2\Im \bar{u}(\nabla-i\A)u$.
This system describes the evolution of a charged nonrelativistic quantum mechanical particle
interacting with the (classical) electro-magnetic field it generates;
$u$ is the wave function of the particle and $(\phi, \A)$ is the electro-magnetic potential.

The solutions of \textbf{MS} have some freedom coming from the gauge invariance,
which is the consequence of the fact that the observables are gauge invariant 
but $u$, $\phi$, and $\A$ themselves are not observables.
Namely, for any function $\lambda: \R^{1+3}\to \R$, \textbf{MS} is invariant under the gauge transform
%$(u',\phi',\A')=(\exp(i\lambda)u, \phi-\partial_t\lambda,\A+\nabla\lambda)$. 
\begin{equation}
(u',\phi',\A')=(\exp(i\lambda)u, \phi-\partial_t\lambda,\A+\nabla\lambda).
\label{gauge}
\end{equation}
By this fact, \eqref{MS1}-\eqref{MS3} itself are not adequate to discuss well-posedness.
For, the uniqueness of the solution clearly does not hold.
To remove this uncertainty, we need to indicate how to choose representative 
elements from  gauge equivalence classes.
Such  conditions are called gauge conditions. 
One of the well-known gauge condition is the Coulomb gauge 
\begin{equation}
\Div \A=0. \label{coulomb}
\end{equation}
In this gauge, \eqref{MS2} and \eqref{MS3} become
\begin{equation}
-\Delta \phi=\rho(u),\quad (\partial_t^2-\Delta) \A=P\J(u,\A),
\label{MS4}
\end{equation}
where $P=\I-\nabla\Div\Delta^{-1}$ is the projection onto the solenoidal subspace.
The first equation in \eqref{MS4} is easily solved by the Newtonian potential.
Therefore \textbf{MS} in the Coulomb gauge (\textbf{MS-C}) is expressed by 
\begin{gather*}
i\partial_tu=(\calh(\A)+\phi(u))u, \quad %\\ %& u(0)=u_0 \\
(\partial_t^2-\Delta) \A =P\J(u,\A),  %& \A(0)=\A_0,\ \ \partial_t\A(0)=\A_1,\ \ \Div \A_0=\Div \A_1=0,
\end{gather*}
where $\phi(u)=(-\Delta)^{-1}|u|^2$, % $u_0$, $\A_0$, $\A_1$ are initial data, 
and the Coulomb gauge condition \eqref{coulomb} is required.
In this gauge $\phi$ does not need the initial datum.
The condition \eqref{coulomb} is conserved if the initial data $\A (0)$ and $\partial_t \A (0)$
satisfy \eqref{coulomb}.
Therefore we consider the time local well-posedness of \textbf{MS-C} with initial data
\begin{equation}
(u(0), \A(0), \partial_t \A(0)) =(u_0,\A_0,\A_1)\in X^{s,\sigma}, \label{eq:IDC}
\end{equation}
where 
$X^{s,\sigma}= \{(u_0,\A_0,\A_1)\in H^s\oplus H^\sigma\oplus H^{\sigma-1} ;
\Div \A_0=\Div \A_1=0\}.$

Our purpose in this paper is to show the local well-posedness for \textbf{MS-C} in the Sobolev space
as wide as possible by the contraction mapping principle.
The main theorem is the following.
%Principally, the most difficult point is to overcome the loss of derivative coming from the term
%$\A\cdot u$.

\begin{theorem}
\label{thmsc}
Let %$s, \sigma$ satisfy 
$s\ge5/3$ and
$\max\{4/3,s-2,(2s-1)/4\} \le \sigma\le \min\{s+1, (5s-2)/3\}$
with $(s,\sigma)\neq (5/2,7/2), (7/2,3/2)${\rm.}
Then for any $(u_0,\A_0,\A_1)\in X^{s,\sigma}${\rm ,}
there exists $T>0$ such that \textbf{\textup{MS-C}} with initial condition \eqref{eq:IDC} 
has a unique solution $(u,\A)$ satisfying
$(u,\A, \partial_t\A)\in C([0,T]; X^{s,\sigma})${\rm .}
Moreover if $\sigma\ge \max\{(s-1),(2s+1)/4\}$ with $(s,\sigma)\neq (5/2,3/2)${\rm ,}
then the map $(u_0,\A_0,\A_1) \mapsto (u,\A,\partial_t\A)$ is continuous 
as a map from $X^{s,\sigma}$ to $C([0,T]; X^{s,\sigma})${\rm .}
\end{theorem}

\vspace{15pt}

%%%%%%%%%%%%%%%%%%%% Figure 
%WinTpicVersion2.15
\unitlength 0.1in
\begin{picture}(43.88,31.22)(19.51,-43.90)
% VECTOR 2 0 3 0
% 2 3126 4499 6333 4499
% 
\special{pn 8}%
\special{pa 3126 4099}%
\special{pa 6333 4099}%
\special{fp}%
\special{sh 1}%
\special{pa 6333 4099}%
\special{pa 6266 4079}%
\special{pa 6280 4099}%
\special{pa 6266 4119}%
\special{pa 6333 4099}%
\special{fp}%
% VECTOR 2 0 3 0
% 2 3423 4790 3423 1874
% 
\special{pn 8}%
\special{pa 3423 4390}%
\special{pa 3423 1474}%
\special{fp}%
\special{sh 1}%
\special{pa 3423 1474}%
\special{pa 3403 1541}%
\special{pa 3423 1527}%
\special{pa 3443 1541}%
\special{pa 3423 1474}%
\special{fp}%
% LINE 2 0 3 0
% 2 4064 1869 4064 4785
% 
\special{pn 8}%
\special{pa 4064 1469}%
\special{pa 4064 4385}%
\special{fp}%
% LINE 2 0 3 0
% 2 3121 3969 6329 3969
% 
\special{pn 8}%
\special{pa 3121 3569}%
\special{pa 6329 3569}%
\special{fp}%
% DOT 2 0 3 0
% 5 3806 4499 4200 4499 4589 4499 4972 4499 5357 4499
% 
\special{pn 8}%
\special{sh 1}%
\special{ar 3806 4099 10 10 0  6.28318530717959E+0000}%
\special{sh 1}%
\special{ar 4200 4099 10 10 0  6.28318530717959E+0000}%
\special{sh 1}%
\special{ar 4589 4099 10 10 0  6.28318530717959E+0000}%
\special{sh 1}%
\special{ar 4972 4099 10 10 0  6.28318530717959E+0000}%
\special{sh 1}%
\special{ar 5357 4099 10 10 0  6.28318530717959E+0000}%
% DOT 2 0 3 0
% 2 5750 4499 6134 4499
% 
\special{pn 8}%
\special{sh 1}%
\special{ar 5750 4099 10 10 0  6.28318530717959E+0000}%
\special{sh 1}%
\special{ar 6134 4099 10 10 0  6.28318530717959E+0000}%
% DOT 2 0 3 0
% 6 3423 4105 3417 3721 3423 3332 3417 2938 3423 2549 3423 2161
% 
\special{pn 8}%
\special{sh 1}%
\special{ar 3423 3705 10 10 0  6.28318530717959E+0000}%
\special{sh 1}%
\special{ar 3417 3321 10 10 0  6.28318530717959E+0000}%
\special{sh 1}%
\special{ar 3423 2932 10 10 0  6.28318530717959E+0000}%
\special{sh 1}%
\special{ar 3417 2538 10 10 0  6.28318530717959E+0000}%
\special{sh 1}%
\special{ar 3423 2149 10 10 0  6.28318530717959E+0000}%
\special{sh 1}%
\special{ar 3423 1761 10 10 0  6.28318530717959E+0000}%
% LINE 2 0 3 0
% 2 3903 4785 6333 2355
% 
\special{pn 8}%
\special{pa 3903 4385}%
\special{pa 6333 1955}%
\special{fp}%
% LINE 2 0 3 0
% 2 3131 4392 5658 1869
% 
\special{pn 8}%
\special{pa 3131 3992}%
\special{pa 5658 1469}%
\special{fp}%
% LINE 2 0 3 0
% 2 4972 2161 3393 4790
% 
\special{pn 8}%
\special{pa 4972 1761}%
\special{pa 3393 4390}%
\special{fp}%
% LINE 2 0 3 0
% 2 4972 2161 5142 1869
% 
\special{pn 8}%
\special{pa 4972 1761}%
\special{pa 5142 1469}%
\special{fp}%
% CIRCLE 2 0 2 0
% 4 4389 3133 4409 3133 4409 3133 4409 3133
% 
\special{pn 8}%
\special{sh 0}%
\special{ar 4389 2733 20 20  0.0000000 6.2831853}%
% LINE 2 0 3 0
% 2 3423 4591 6329 3138
% 
\special{pn 8}%
\special{pa 3423 4191}%
\special{pa 6329 2738}%
\special{fp}%
% LINE 2 0 3 0
% 2 3423 4591 3126 4736
% 
\special{pn 8}%
\special{pa 3423 4191}%
\special{pa 3126 4336}%
\special{fp}%
% CIRCLE 2 0 2 0
% 4 4778 3910 4798 3906 4798 3906 4798 3906
% 
\special{pn 8}%
\special{sh 0}%
\special{ar 4778 3510 20 20  0.0000000 6.2831853}%
% LINE 2 2 3 0
% 2 6339 1966 3514 4785
% 
\special{pn 8}%
\special{pa 6339 1566}%
\special{pa 3514 4385}%
\special{dt 0.045}%
\special{pa 3514 4385}%
\special{pa 3514 4385}%
\special{dt 0.045}%
% LINE 2 2 3 0
% 2 6339 2938 3131 4537
% 
\special{pn 8}%
\special{pa 6339 2538}%
\special{pa 3131 4137}%
\special{dt 0.045}%
\special{pa 3131 4137}%
\special{pa 3132 4137}%
\special{dt 0.045}%
% CIRCLE 2 0 2 0
% 4 4389 3910 4404 3910 4404 3910 4404 3910
% 
\special{pn 8}%
\special{sh 0}%
\special{ar 4389 3510 20 20  0.0000000 6.2831853}%
% STR 2 0 3 0
% 3 3427 1705 3427 1753 5 0
% $\sigma$
\put(34.2700,-13.5300){\makebox(0,0){$\sigma$}}%
% STR 2 0 3 0
% 3 6411 4455 6411 4503 5 0
% $s$
\put(64.1100,-41.0300){\makebox(0,0){$s$}}%
% STR 2 0 3 0
% 3 4200 4542 4200 4591 5 0
% $2$
\put(42.0000,-41.9100){\makebox(0,0){$2$}}%
% STR 2 0 3 0
% 3 3806 4547 3806 4596 5 0
% $1$
\put(38.0600,-41.9600){\makebox(0,0){$1$}}%
% STR 2 0 3 0
% 3 6577 3055 6577 3103 5 0
% $\sigma=\frac{2s-1}{4}$
\put(67.5700,-27.0300){\makebox(0,0){$\sigma=\frac{2s-1}{4}$}}%
% STR 2 0 3 0
% 3 6127 2875 6127 2923 2 0
% $\sigma=\frac{2s+1}{4}$
\put(64.4700,-25.2300){\makebox(0,0)[lb]{$\sigma=\frac{2s+1}{4}$}}%
% STR 2 0 3 0
% 3 6397 2245 6397 2293 5 0
% $\sigma=s-2$
\put(67.9700,-18.9300){\makebox(0,0){$\sigma=s-2$}}%
% STR 2 0 3 0
% 3 6397 1795 6397 1843 5 0
% $\sigma=s-1$
\put(67.9700,-14.4300){\makebox(0,0){$\sigma=s-1$}}%
% STR 2 0 3 0
% 3 5677 1705 5677 1753 5 0
% $\sigma=s+1$
\put(56.7700,-13.5300){\makebox(0,0){$\sigma=s+1$}}%
% STR 2 0 3 0
% 3 5137 1705 5137 1753 5 0
% $\sigma=\frac{5s-2}{3}$
\put(49.3700,-13.5300){\makebox(0,0){$\sigma=\frac{5s-2}{3}$}}%
% STR 2 0 3 0
% 3 4057 1705 4057 1753 5 0
% $s=\frac{5}{3}$
\put(40.5700,-13.5300){\makebox(0,0){$s=\frac{5}{3}$}}%
% STR 2 0 3 0
% 3 6397 3865 6397 3913 5 0
% $\sigma=4/3$
\put(66.9700,-35.5300){\makebox(0,0){$\sigma=4/3$}}%
% STR 2 0 3 0
% 3 4584 4542 4584 4591 5 0
% $3$
\put(45.8400,-41.9100){\makebox(0,0){$3$}}%
% STR 2 0 3 0
% 3 4972 4542 4972 4591 5 0
% $4$
\put(49.7200,-41.9100){\makebox(0,0){$4$}}%
% STR 2 0 3 0
% 3 5357 4542 5357 4591 5 0
% $5$
\put(53.5700,-41.9100){\makebox(0,0){$5$}}%
% STR 2 0 3 0
% 3 5745 4542 5745 4591 5 0
% $6$
\put(57.4500,-41.9100){\makebox(0,0){$6$}}%
% STR 2 0 3 0
% 3 6134 4542 6134 4591 5 0
% $7$
\put(61.3400,-41.9100){\makebox(0,0){$7$}}%
% STR 2 0 3 0
% 3 3325 4537 3325 4586 5 0
% $0$
\put(32.5500,-41.8600){\makebox(0,0){$0$}}%
% STR 2 0 3 0
% 3 3315 4056 3315 4105 5 0
% $1$
\put(33.1500,-37.0500){\makebox(0,0){$1$}}%
% STR 2 0 3 0
% 3 3320 3667 3320 3716 5 0
% $2$
\put(33.2000,-33.1600){\makebox(0,0){$2$}}%
% STR 2 0 3 0
% 3 3320 3284 3320 3332 5 0
% $3$
\put(33.2000,-29.3200){\makebox(0,0){$3$}}%
% STR 2 0 3 0
% 3 3320 2890 3320 2938 5 0
% $4$
\put(33.2000,-25.3800){\makebox(0,0){$4$}}%
% STR 2 0 3 0
% 3 3320 2501 3320 2549 5 0
% $5$
\put(33.2000,-21.4900){\makebox(0,0){$5$}}%
% STR 2 0 3 0
% 3 3320 2112 3320 2161 5 0
% $6$
\put(33.2000,-17.6100){\makebox(0,0){$6$}}%
% STR 2 0 3 0
% 3 4381 3829 4381 3877 3 0
% $(\frac{5}{2},\frac{3}{2})$
\put(44.8100,-34.4700){\makebox(0,0)[rb]{$(\frac{5}{2},\frac{3}{2})$}}%
% STR 2 0 3 0
% 3 4597 3145 4597 3193 5 0
% $(\frac{7}{2},\frac{5}{2})$
\put(46.3700,-27.9300){\makebox(0,0){$(\frac{5}{2},\frac{7}{2})$}}%
% STR 2 0 3 0
% 3 4633 3833 4633 3882 5 0
% $(\frac{7}{2},\frac{3}{2})$
\put(49.5300,-36.8200){\makebox(0,0){$(\frac{7}{2},\frac{3}{2})$}}%
\put(47.8400,-46.9100){\makebox(0,0){Figure 1}}%
\end{picture}%
%%%%%%%%%%%%%%%%%%% end Figure

\vspace{15pt}

\newpage
\begin{remark}
(1)
$T$ depends only on $s,\sigma$ and $\|(u_0,\A_0,\A_1);X^{s,\sigma}\|$.

(2)
For any $s$ and $\sigma$ satisfying the assumption above for the unique existence of the solution, 
the map $(u_0,\A_0,\A_1) \mapsto (u,\A,\partial_t\A)$ 
is continuous in weak-star sense.
\end{remark}

Nakamitsu-Tsutsumi \cite{NT} showed the time local well-posedness in $X^{s,s}$ with $s>5/2$.
In fact, they treated the case of Lorentz gauge mentioned below, 
but the Coulomb gauge case can be treated analogously.
Generally, the most difficult point of the treatment of \textbf{MS-C} is 
to overcome the loss of derivative
which may be caused by the term $\A\nabla u$ in \eqref{MS1}.
In \cite{NT} it is done by usual energy method.
The fact that $\Re\langle \A\nabla u,u\rangle_{H^s}=-\int |\Omega^s u|^2\Div\A dx=0$,
where $\Omega=(1-\Delta)^{1/2}$, is used to obtain a differential inequality 
$d\|u;H^s\|/dt \lesssim \| \A;H^s\| \| \partial u\|_\infty +\| \partial \A \|_\infty \| u;H^s\|$.
The assumption $s>5/2$ is needed to treat $\| \partial u \|_\infty$ and $\| \partial \A \|_\infty$.
In order to refine the result, in the present paper 
we  derive the estimate for $\| \calh u ;H^{s-2}\|$ or $\| \partial_t u;H^{s-2} \|$
instead of $\|u;H^s\|$ itself. % since \eqref{MS1} is first order in $t$ and second order in $x$.
Then the self-adjointness of $\calh$ in $L^2$ helps us to overcome the loss of derivative.
We also remark that the energy inequality for the wave equation in $H^\sigma$ requires that the inhomogeneous 
term belongs to $H^{\sigma-1}$, from which the assumption $\sigma\le s$ seems to be needed.
However, we actually need a weaker condition for $\sigma$ by the use of the projection $P$
(see Lemma \ref{wavest}).
On the other hand, Guo-Nakamitsu-Strauss \cite{GNS} constructed a time global solution in $X^{1,1}$ 
although they did not show the uniqueness.
Indeed, \textbf{MS-C} has the conservation laws of the charge and the energy from which we can obtain the boundedness of 
$\| (u,\A,\partial_t \A);X^{1,1}\|$. Therefore this result is obtained by the compactness method.
Our result fills some part of the gap between \cite{GNS} and \cite{NT} but not completely.

Next we consider the Lorentz gauge 
\begin{equation}
\partial_t\phi+\Div \A=0.
\label{lorentz}
\end{equation}
\textbf{MS} in the Lorentz gauge (\textbf{MS-L}) is expressed as
\begin{gather*}
i\partial_tu=(\calh(\A)+\phi)u, \quad %\\
(\partial_t^2-\Delta) \phi =\rho(u), \quad %\\
(\partial_t^2-\Delta) \A =\J(u,\A).
\end{gather*}
In this case, we need the initial data
\begin{equation}\label{eq:IDL}
(u(0),\phi(0),\partial_t \phi(0), \A(0), \partial_t \A(0))
=(u_0,\phi_0,\phi_1,\A_0,\A_1)\in Y^{s,\sigma}.
\end{equation}
Here
\begin{align*}
Y^{s,\sigma}&=\{ (u_0,\phi_0,\phi_1,\A_0,\A_1) \in 
H^s\oplus H^\sigma\oplus H^{\sigma-1}\oplus H^\sigma\oplus H^{\sigma-1}; \\
&\quad \Div \A_0+\phi_1=\Div \A_1+\Delta \phi_0+|u_0|^2=0 \}. 
\end{align*}
The condition \eqref{lorentz} is conserved if the initial datum belongs to $Y^{s,\sigma}$.
The result for \textbf{MS-L} is the following.

\begin{theorem}
\label{thmsl}
Let $s\ge5/3$ and
$\max\{4/3,s-1 \} \le \sigma\le \min\{s+1, (5s-2)/3\}$
with $(s,\sigma)\neq (5/2,7/2)${\rm .}
Then for any
$(u_0,\phi_0,\phi_1,\A_0,\A_1)\in Y^{s,\sigma}${\rm ,}
there exists $T>0$ such that 
{\bf MS-L} with initial condition \eqref{eq:IDL} has a unique solution 
$(u,\phi,\A)$ satisfying 
\[(u,\phi,\partial_t \phi,\A,\partial_t \A) \in C([0,T];Y^{s,\sigma}). \] %{\rm .}
Moreover{\rm ,} if $\sigma\ge(2s+1)/4$ with $(s,\sigma)\neq (5/2,3/2)${\rm ,} 
then the map
$(u_0,\phi_0,\phi_1,\A_0,\A_1)\mapsto (u,\phi,\partial_t \phi,\A,\partial_t \A)$ 
is continuous as a map from $Y^{s,\sigma}$ to $C([0,T]; Y^{s,\sigma})${\rm .}
\end{theorem}

We can also treat the temporal gauge, namely 
\begin{equation}\label{temporal}
\phi=0.
\end{equation}
In this gauge \textbf{MS} becomes the following system, which is referred to as \textbf{MS-T}:
\begin{gather*}
i\partial_tu=\calh(\A)u, \quad
(\partial_t^2-\Delta) \A +\nabla \Div \A=\J(u,\A).
\end{gather*}
For \textbf{MS-T}, we need the initial data
\begin{equation}\label{eq:IDT}
(u(0),\A(0), \partial_t \A(0))=(u_0,\A_0,\A_1)\in \tilde{Y}^{s,\sigma},
\end{equation}
where $\tilde{Y}^{s,\sigma}=\{(u_0,\A_0,\A_1)\in H^s\oplus H^\sigma \oplus H^{\sigma-1};-\Div \A_1=|u_0|^2\}$.
\begin{theorem}
\label{thmst}
Let $s\ge5/3$ and
$\max\{4/3,s-1 \} \le \sigma\le \min\{s+1, (5s-2)/3\}$
with $(s,\sigma)\neq (5/2,7/2)${\rm .}
Then there exists $T>0$ such that  
{\bf MS-T} with initial condition \eqref{eq:IDT} has a unique solution 
$(u,\A)$ satisfying 
$(u,\A,\partial_t \A) \in C([0,T];\tilde{Y}^{s,\sigma})${\rm .}
Moreover{\rm ,} if $\sigma\ge(2s+1)/4$ with $(s,\sigma)\neq (5/2,3/2)${\rm ,} 
then the map
$(u_0,\A_0,\A_1)\mapsto (u,\A,\partial_t \A)$
is continuous as a map from $\tilde{Y}^{s,\sigma}$ to $C([0,T]; \tilde{Y}^{s,\sigma})${\rm .}\end{theorem}

This paper is organized as follows.
In Section 2, we introduce some elementary estimates required in this paper. 
In Section 3, we construct the evolution operator for the linear Schr\"odinger equation.
In Section 4, we prepare a priori estimates for the solutions of linearized equation.
In Sections 5 and 6, we prove Theorem \ref{thmsc}  by the contraction mapping principle 
except for the continuous dependence of the solutions on the data, 
which is proved in Section 7. 
In Section 8, We prove Theorems \ref{thmsl} and \ref{thmst}.

We conclude this section by giving the notation used in this paper.
$\omega = (-\Delta)^{1/2}$ and $\Omega =(1-\Delta)^{1/2}$.
$L^p=L^p(\R^3)$ is the usual Lebesgue space 
and its norm is denoted by $\|\cdot\|_p$ or $\|\,\cdot\,;1/p\|$.
$p'=p/(p-1)$ is the dual exponent of $p$. 
This symbol is used only for Lebesgue exponents.
$H^{s,p}=\{ \phi \in \mathcal{S}'(\R^3); \| \Omega^s \phi \|_p <\infty \}$ 
is the usual Sobolev space.
For any interval $I \subset \R$ and Banach space $X$, $L^p(I;X)$ denotes the space of 
$X$-valued strongly measurable functions on $I$ whose $X$-norm belong to $L^p(I)$.
This space is often abbreviated to $L^pX$ when we fix the time interval $I$.
$W^{m,p}(I;X)$ denotes the space of functions in $L^p(I;X)$ 
whose derivatives up to the $(m-1)$-times are locally absolutely continuous 
and the derivatives up to the $m$-times belong to $L^p(I;X)$.
%
%$W^{1,p}(I;X)$ denotes the space of $X$-valued absolutely continuous functions on $I$ 
%whose derivative belong to $L^p(I;X)$.
%$W^{m,p}(I;X)=\{ u(t) ; \partial_t^j u(t) \in W^{1,p}(I;X) \text{ for } 0\le j \le m-1 \}$ 
%
The inequality $a\lesssim b$ means $a\le Cb$, 
where $C$ is a positive constant that is not essential.
%For any real number $a$, $[a]$ denotes the integral part of $a$.
$\langle a\rangle=\sqrt{1+a^2}$. 
$a\vee b$ and $a\wedge b$ denote the maximum and the minimum of $a$ and $b$ respectively.
We use the following unusual but convenient symbol:
$a_+$ means $a\vee0$ if $a\neq0$, whereas $0_+$ means a sufficiently small positive number.
Namely
$b\ge a_+$ means $b\ge a\vee 0$ if $a\neq 0$ and $b>0$ if $a=0$.
It is useful to express 
sufficient conditions for Sobolev type embeddings $H^{s,r} \embedding L^p$
by the inequality $(1/r-s/3)_+ \le 1/p\le 1/r$.

%This symbol may be unusual but it is convinient to express the following condition:
%Certainly $0_+$ is not uniquely determined, 
%but there is no confusion since it is only used to express 
%a sufficient condition for Sobolev type embeddings $H^{s,r} \embedding L^p$
%by the inequality $(1/r-s/3)_+ \le 1/p\le 1/r$.
%We consider \textbf{MS-C} and we put $\phi(u)=\omega^{-2}|u|^2$ 
%in the following except for the section 8, 
%where $\omega$ denotes $\sqrt{-\Delta}$.

%%%%%%%%%%%%%%%%% End Introduction  %%%%%%%%%%%%%%%%%%%%%%%%%%%%%%%%%%

%%%%%%%%%%%%%%%%% Preliminaries   %%%%%%%%%%%%%%%%%%%%%%%%%%%%%%%%%%%%
\section{Preliminaries}
\begin{lemma} 
\label{phies}
Let $s, s_1, s_2, s_3$ satisfy $0\le s\le s_3${\rm ,} 
$s_1\wedge s_2\ge(s-2)\vee0${\rm ,} and $s_1+s_2>0${\rm .}
Let 
\[
s_1+s_2+s_3\wedge (3/2)\ge s+1
\]
and the inequality be strict if {\rm (1)} $s_j=3/2$ for some $1\le j\le 3$
or {\rm (2)} $s=s_3<3/2${\rm .}
Then the following estimate holds{\rm:}
\begin{equation}\label{eq:Hartreenlest}
\| \omega^{-2} (u_1 u_2) u_3 ;H^s \| 
\lesssim 
\prod_{j=1}^3 \| u_j ; H^{s_j} \|. 
\end{equation}
\end{lemma}

\Proof
By Leibniz's rule the left-hand side of \eqref{eq:Hartreenlest} 
is bounded by some constant times 
\begin{equation}\label{eq:HLeibniz}
\| \omega^{-2} (u_1 u_2) \|_{p_1} \| u_3 ;H^{s,p_2} \|
+
\| \omega^{-2+s} (u_1 u_2) \|_{p_3} \| u_3 \|_{p_4}
\equiv\textrm{I}+\textrm{II}
\end{equation}
with $1/2=1/p_1+1/p_2=1/p_3+1/p_4$.
We begin with the treatment of the first term.
We choose $p_2$ as large as possible provided that % on conditions that
$H^{s_3} \embedding H^{s,p_2}$ and that the operator 
$\omega^{-2}:L^\nu \to L^{p_1}$
($0<1/\nu=1/p_1+2/3<1$) is bounded by virtue of 
the Hardy-Littlewood-Sobolev inequality.
Certainly if $p_1 =\infty$, we use H\"older's inequality instead.
With such a choice of $(p_1,p_2)$ we have
$\textrm{I} \lesssim \| u_1 u_2 \|_\nu \| u_3 ; H^{s_3}\|$.
Next we apply H\"older's inequality  and Sobolev's inequality
for the first factor in order to obtain
$\| u_1 u_2 \|_\nu \lesssim \prod_{j=1}^2 \| u_j ; H^{s_j}\|$.
Then we obtain
\[ \textrm{I}\lesssim \prod_{j=1}^3 \| u_j ; H^{s_j}\|. \]
Specifically, if $s_3=s$, we choose $p_1=\infty$, $p_2 =2$. 
Then we need $s_1+s_2>1$ to obtain the estimate above.
If $0<s_3-s<1$, we choose %$1/p_1=(s_3-s)/3$, 
$1/p_2 =1/2 -(s_3-s)/3$,
and then we need $s_1+s_2+s_3\ge s+1$.
If $s_3-s\ge 1$, we choose $p_1 =3+0$
so that $\nu=1+0$, and then we need $s_1+s_2>0$.

We proceed to the treatment of the second term,
which is essentially similar to that of the first term,
but we have to divide the proof into cases in a different fashion.
If $s=0$, we do not need the estimate for $\textrm{II}$.
If $0<s\le (s_3\wedge 3/2)-1$, 
we choose $1/p_3 =(1+s)/3 -0$, $1/p_4 = 1/2 -(1+s)/3+0$.
Then we have 
$\| \omega^{-2+s} (u_1 u_2) \|_{p_3} 
\lesssim \prod_{j=1}^2 \| u_j; H^{s_j} \|$ 
provided $s_1+s_2>0$.
Since we also have $\| u_3 \|_{p_4} \lesssim \| u_3; H^{s_3} \|$,
we obtain 
\begin{equation}\label{eq:IIest}
\textrm{II}\lesssim \prod_{j=1}^3 \| u_j ; H^{s_j}\|.
\end{equation}
If $(s_3\wedge 3/2)-1 <s \le 2$ and if $s_3 \neq 3/2$,
we choose $1/p_3 = (s_3/3)\wedge (1/2)$,
$1/p_4 = (1/2 -s_3/3) \vee 0$. 
Then we obtain \eqref{eq:IIest}
provided $s_1+s_2+(s_3\wedge 3/2)\ge s+1$.
We can similarly estimate $\textrm{II}$ even if $s_3 = 3/2$,
but the limiting case $s_1+s_2+s_3=s-1/2$ is excluded from the 
sufficient condition for the estimate
because of the exception of Sobolev's embedding theorem, 
namely $H^{3/2} \not\subset L^{\infty}$.
If $s>2$, we choose $p_3 = 2, p_4 =\infty$.
Then we obtain \eqref{eq:IIest} by virtue of
Leibniz's formula together with Sobolev's inequality.
The exceptionally prohibited case in the statement of the lemma 
comes from the exception of Sobolev's embedding theorem. 
\ebox

\begin{lemma}
\label{nablaes}
{\rm (1)} Let $\sigma\ge s\vee(1/2)\vee(-s-1)$ and $(s,\sigma)\neq (1/2,1/2), (-3/2,1/2)${\rm .}
Then 
\begin{equation}
\|(\nabla-i\A)v;H^s\|\lesssim \|v;H^{s+1}\|\langle\|\A;H^\sigma\|\rangle
\label{nablaes1}
\end{equation}
for any $v$ and $\A${\rm .}

{\rm (2)} 
Let $s, s_1, s_2$ satisfy $0\le s\le s_1\wedge s_2${\rm ,} $s_1+s_2\ge s+3/2${\rm ,}
$(s_1,s_2)\neq(s,3/2),(3/2,s)${\rm ,} and let $\sigma$ satisfy
$\sigma\ge s_2\vee (1/2)${\rm ,} $(\sigma,s_2)\neq (1/2,1/2)${\rm .}
Then
\begin{equation}
\|w(\nabla-i\A)v;H^s\| \lesssim \|w;H^{s_1}\|\|v;H^{s_2+1}\|\langle\|\A;H^\sigma\|\rangle 
\label{nablaes2}
\end{equation}
for any $w, \A$ and $v${\rm .}
\end{lemma}

\Proof 
(1) If $s\ge0$, by the Leibniz formula and the Sobolev inequality, 
we have $\|\A v;H^s\|\lesssim \|\A;H^\sigma\| \|v;H^{s+1}\|$ 
with $\sigma\ge(1/2)\vee s, (s,\sigma)\neq(1/2,1/2)$.
If $-1\le s<0$, again by the Sobolev inequality
\[
\|\A v;H^s\| \lesssim \|\A\|_3\|v;1/2-(s+1)/3\| \lesssim \|\A;H^\sigma\| \|v;H^{s+1}\|
\]
with $\sigma\ge1/2$.
If $s<-1$, we use duality.
We have by the result for $s\ge0$
\begin{align*}
|\langle\A v,\psi\rangle| & \le \|v;H^{s+1}\|\|\A\psi;H^{-s-1}\| \\
&\lesssim \|v;H^{s+1}\|\|\A;H^\sigma\|\|\psi;H^{-s}\|.
\end{align*}
This estimates yields $\|\A v;H^s\|\lesssim \|\A;H^\sigma\| \|v;H^{s+1}\|$.
Consequently we have this inequality for all $s,\sigma$.
Using this estimate with the trivial estimate $\|\nabla v;H^s\|\lesssim \|v;H^{s+1}\|$, 
we obtain \eqref{nablaes1}.

(2) By the Leibniz formula, we have
\[
\|w(\nabla-i\A)v;H^s\| \lesssim \|w\|_{r_1}\|(\nabla-i\A)v;H^{s,2r_1/(r_1-2)}\|
+\|w;H^{s,2r_2/(r_2-2)}\|\|(\nabla-i\A)v\|_{r_2},
\]
where $1/r_j=(1/2-s_j/3)_+$, $j=1,2$.
Under the assumption for $s,s_1$ and $s_2$, the right-hand side does not exceed
some positive constant times $\|w;H^{s_1}\|\|(\nabla-i\A)v;H^{s_2}\|$ 
by virtue of the Sobolev inequality.
Therefore we obtain \eqref{nablaes2} by (1).
\ebox
\thmskip

We define $\Gamma$ by 
\[
\Gamma \equiv 
\left\{(s,\sigma);
\begin{array}{l}
 s\ge0, \sigma\ge (3/4-s/2)\vee(1/2)\vee (s/2-1/4)\vee(s-2), \\
(s,\sigma)\neq (7/2,3/2),(3/2,1/2),(1/2,1/2)
\end{array}
\right\}.
\]

\begin{lemma}
\label{equinorm}
Let $(s,\sigma)\in \Gamma${\rm ,} $\Div \A=0${\rm .}
Then 
\begin{equation}
\|(\calh(\A)+\phi(u))v;H^{s-2}\| 
\lesssim \|v;H^s\|\langle\|\A;H^\sigma\|\vee \|u;H^{(s-1) \vee0}\|\rangle^2.
\label{Hphi}
\end{equation}
Moreover if $s>0${\rm ,} $\sigma>(1/2) \vee(3/4-s/2)${\rm ,} then 
\begin{equation}
\|v;H^s\| \lesssim\|(\calh(\A)+\phi(u))v;H^{s-2}\|+\langle\|\A;H^\sigma\| \vee \|u;H^{(s-1) \vee0}\| \rangle^\alpha \|v\|_2,
\label{equi2}
\end{equation}
where $\alpha=\alpha(s,\sigma)$ is a positive constant independent of $v$ and $\A${\rm .}
\end{lemma}

\Proof
First we show 
\begin{equation}
\|(2i\A\nabla +\A^2)v;H^{s-2}\| \lesssim \|v;H^s\|\langle\|\A;H^\sigma\|\rangle^2.
\label{equi1}
\end{equation}
%To prove this, we should estimate $\|\A\nabla v;H^{s-2}\|$ and $\|\A^2v;H^{s-2}\|$.
%We have 
%\[
%\|(\nabla-i\A)^2v;H^{s-2}\| \lesssim \|\Delta v;H^{s-2}\|+\|\A\nabla v;H^{s-2}\|+\|\A^2v;H^{s-2}\|.
%\]
For $0\le s\le 2$, by the Sobolev inequality, 
\begin{align*}
\|\A \nabla v;H^{s-2}\| &\lesssim
\begin{cases}
\|\A \nabla v; 1/2-(s-2)/3\| & \text{if } 1\le s\le 2, \\
\|\A v; H^{s-1}\| & \text{if } 0\le s< 1
\end{cases} \\
&\lesssim \|\A\|_3\|v;H^s\|.
\end{align*}
Here we have used $\Div \A =0$ when $0\le s<1$.
We also have
\begin{align*}
\|\A^2v;H^{s-2}\| \lesssim \|\A^2v \|_r 
&\lesssim \|\A;(1/2-\sigma/3)_+\|^2 \|v;1/r-2(1/2-\sigma/3)_+\| \\
& \lesssim \|\A;H^\sigma\|^2\|v;H^s\|,
\end{align*}
where $1/2\le 1/r \le 1-(1/2-(2-s)/3)_+$, $(1/2-s/3)_+\le 1/r-2(1/2-\sigma/3)_+\le 1/2$.
Such $r$ exists if $(s,\sigma) \in \Gamma$ with $0\le s\le 2$. 
These estimates imply \eqref{equi1} for $0\le s\le 2$. 
For $s>2$, by the Leibniz formula and the Sobolev inequality, we have 
$\|\A\nabla v;H^{s-2}\| \lesssim \|\A;H^\sigma\| \|v;H^s\|$ and 
\begin{align*}
\|\A^2v;H^{s-2}\| &\lesssim \|\A;H^{s-2,2r/(r-2)}\| \|\A\|_{r}\|v\|_\infty
+\|\A\|_{r}^2\|v;H^{s-2,2r/(r-4)}\| \\
&\lesssim \|\A;H^\sigma\|^2\|v;H^s\|,
\end{align*}
where $\sigma\ge s-2+3/r$, $1/r=(1/2-\sigma/3)_+$.
Thus \eqref{equi1} has been established.
We remark that actually we have 
\begin{align}
\|\A\nabla v;H^{s-2}\| &\lesssim \|\A;H^\sigma\|\|v;H^{s-\delta}\|, \label{eq:sharpADv}\\
\|\A^2v;H^{s-2}\|&\lesssim \|\A;H^\sigma\|^2\|v;H^{s-\delta}\|, \label{eq:sharpAAv}
\end{align}
if $s>0$, $\sigma>1/2$ and $\sigma>(3/4-s/2)$.
Here $\delta$ is a sufficiently small positive number. 
These inequalities will be used to prove \eqref{equi2}.
The estimates \eqref{Hphi} follows from \eqref{equi1} and the inequality
\begin{equation}\label{eq:phies1}
\|\phi(u)v;H^{s-2}\|\lesssim \|u;H^{(s-1)\vee 0}\|^2\|v;H^{(s-1)_+}\|.
\end{equation}
If $s\ge2$, this inequality follows from Lemma \ref{phies} directly. If $s<2$, this follows from the duality argument such as 
\[
|\langle\phi(u)v, \psi\rangle| = |\langle u, \omega^{-2}(\bar{v}\psi)u\rangle| 
%\le \|u\|_2 \| \omega^{-2}(\bar{v}\psi)u\|_2 
\lesssim \|u\|_2\|u;H^{(s-1)\vee 0}\|\|v;H^{(s-1)_+}\|\|\psi;H^{2-s}\|.
\]
Next we show \eqref{equi2}. 
Clearly
\[ \|v;H^s\| \lesssim \|(\calh(\A)+\phi(u))v;H^{s-2}\|+\|(2i\A\nabla+\A^2+\phi(u))v;H^{s-2}\|+\| u \|_2. \]
We apply the interpolation inequality
\begin{equation}
a\|v;H^{s-\delta}\|\lesssim a\|v\|_2^{\delta/s}\|v;H^s\|^{(s-\delta)/s}
\lesssim \varepsilon\|v;H^s\|+C(\varepsilon) a^{s/\delta}\|v\|_2
\label{interpolation}
\end{equation}
to \eqref{eq:sharpADv}-\eqref{eq:phies1},
where $a>0$ is a constant. Then we obtain \eqref{equi2} by taking $\varepsilon$ sufficiently small. 
\ebox

\begin{remark}
%(1)
Let the assumption for \eqref{equi2} be satisfied and
let $\A (t) \in C(I;H^\sigma)$, $u(t),v(t) \in C(I;H^{s-\delta})$ for some $\delta>0$.
Then by estimates similar to \eqref{eq:sharpADv}-\eqref{eq:phies1}, we have
$\A\nabla v, \A^2 v$, $\phi(u)v \in C(I;H^s)$.
This fact will be used later.
%
%(2) 
%By the proof, we have the inequality
%$\|\A(\nabla-i\B)v;H^{s-2}\| \lesssim \| \A;H^\sigma\| \langle\| \B;H^\sigma\|\rangle \| v;H^s\|$
%under the assumption of the lemma.
\end{remark}

\section{Estimates for solutions to linear Schr\"odinger equations}

For our treatment of \textbf{MS},
we need energy estimates for linear Schr\"odinger equations with
electro-magnetic potentials.
Let $I\subset \R$ be a compact interval,
$\phi:I\times \R^3 \to \R$, $\A :I\times\R^3 \to \R^3$, 
and  $t_0 \in I$.
We consider the equation 
\begin{equation}\label{eq:LS}
i\partial_t v = \mathcal{H}(\A ) v +\phi v
\end{equation}
%where $\mathcal{H}(\A )=-(\nabla -i\A )^2$, 
with initial data
\begin{equation}\label{eq:LSinitial}
v(t_0) =v_0.
\end{equation}
For a while we regard $u$ and $\A $ as known functions,
and consider the linear Cauchy problem 
\eqref{eq:LS}-\eqref{eq:LSinitial}.

Before we proceed to energy estimates, we clarify 
the concept of the solution.
Let $s\ge 0$. A function $v$ is called an $H^s$-solution to \eqref{eq:LS}
in $I$ if $v\in C(I; H^s) \cap W^{1,1}( I; H^{s-2})$ and satisfies
\eqref{eq:LS} almost every $t\in I$.
Moreover, if $v_0 \in H^s$ and $v$ satisfies \eqref{eq:LSinitial},
then $v$ is called an $H^s$-solution to \eqref{eq:LS}-\eqref{eq:LSinitial}.

The $L^2$-norm of $H^1$-solutions are clearly conserved.
For, if $v(t)$ is an $H^1$-solution, 
then $\| v(t) \|_2^2$ is absolutely continuous, and \eqref{eq:LS} yields
\[ d \| v(t) \|_2^2/dt = 2\textrm{Re}
\langle \partial_t v, v \rangle_{H^{-1}\times H^1} =0 \]
for almost every $t\in I$.
This also implies the uniqueness of the $H^1$-solution to 
\eqref{eq:LS}-\eqref{eq:LSinitial}.

For a while we assume the following:
\begin{assumption}\label{A1}\rm
(1) $\A  \in L^\infty (I;H^1)\cap W^{1,1}(I; L^3)$
with $\Div \A =0$;

(2) $\phi =\phi (u) = \omega^{-2} |u|^2$ with 
$u\in L^\infty (I; H^{3/4})$.
\end{assumption}
\begin{lemma}\label{lem:LSH2}
We assume {\rm \ref{A1}.}
Then \eqref{eq:LS}-\eqref{eq:LSinitial} has a unique $H^2$-solution $v${\rm ,}
which satisfies the estimate
\begin{equation}\label{eq:H2estLS}
\begin{split}
&\| \mathcal{H}(\A ) v(t) \|_2 + \langle l \rangle^4 \| v(t) \|_2 \\
&\quad\le \bigl\{ \| \mathcal{H}(\A (t_0)) v_0 \|_2 
  + \langle l \rangle^4 \| v_0 \|_2 \bigr\} \\
&\quad\quad \times \exp \bigl\{ C\| u;L^2(I;H^{3/4}) \|^2 
  +C\| \partial_t \A  ; L^1(I;L^3)\| \bigr\}
\end{split}
\end{equation}
for any $t\in I${\rm .} Here $l=\|\A ;L^\infty (I;H^1)\|${\rm .}
Moreover{\rm ,} if $u\in C(I;L^2)$ then $v\in C^1(I;L^2)${\rm .} 
\end{lemma}
\Proof
We  first prove \eqref{eq:H2estLS} rather formally.
By direct computation, self-adjointness of $\mathcal{H}(\A )$ 
and Schwarz's inequality
\begin{align*}
\frac{1}{2} \frac{d}{dt} \| \mathcal{H} v \|_2^2
&= \textrm{Im} 
\langle \mathcal{H}(\mathcal{H}+\phi)v
+2 \partial_t \A (\nabla -i\A )v, \mathcal{H} v \rangle
\\
&\le\bigl\{ \| \mathcal{H} \phi v \|_2 
+ 2 \| \partial_t \A (\nabla -i\A )v \|_2 \bigr\}
  \| \mathcal{H} v \|_2.
\end{align*}
The quantities in the brackets are estimated 
by Sobolev's inequality, Lemmas \ref{phies} and \eqref{equi2} with $u=0$.
Indeed we have
\begin{align*}
\| \mathcal{H} \phi v \|_2 
&\lesssim \| \phi v;H^2 \| +\langle l \rangle^4 \| \phi v \|_2 \\
&\lesssim \| u;H^{3/4} \|^2 \{ \| v;H^2 \| + \langle l \rangle^4 \| v \|_2 \} \\
&\lesssim \| u;H^{3/4} \|^2 
   \{ \| \mathcal{H} v \|_2 + \langle l \rangle^4 \| v \|_2 \},
\end{align*}
and
\[
\| \partial_t \A (\nabla -i\A )v \|_2
\le \| \partial_t \A  \|_3 \| (\nabla -i\A )v \|_6 
\lesssim \| \partial_t \A  \|_3 
   \{ \| \mathcal{H} v \|_2 + \langle l \rangle^4 \| v \|_2 \}.
\]
We remark that we can take $\alpha=4$ in \eqref{equi2}.
In the last inequality we have used the estimate 
\[
\| \A v \|_6 \le \| \A \|_6 \| v\|_\infty 
\lesssim \| \A ; H^1 \| \| v \|_2^{1/4} \| v ;H^2 \|^{3/4}
\lesssim \| v ;H^2 \| + l^4 \| v \|_2.
\]
Therefore we obtain the differential inequality
\[ \frac{d}{dt} \{ \| \mathcal{H} v \|_2 + \langle l \rangle^4 \| v \|_2 \}
\lesssim \{ \| u;H^{3/4} \|^2 + \| \partial_t \A  \|_3 \}
\{ \| \mathcal{H} v \|_2 + \langle l \rangle^4 \| v \|_2 \},   
\]
which yields \eqref{eq:H2estLS} by virtue of Gronwall's inequality. 
Here we have used the $L^2$-norm conservation law in the light-hand side.
We next prove the existence of the solution.
If $\A $ and $u$ are sufficiently smooth, 
the existence of the solution is proved by Kato's abstract
method~\cite{K70,K73},
or the parabolic regularization technique.
Indeed, the condition $u, \A  \in L^2 (I; H^{s+1/2 +0})$ with $s\ge 2$
will suffice to prove the $H^s$-wellposedness.
To construct the solution under the assumption \ref{A1}, we put $\A _k = \eta_k \ast \A $, $u_k = \eta_k \ast u$,
$k=1,2,\cdots$,
and consider the problem \eqref{eq:LS}-\eqref{eq:LSinitial}
with $(u, \A )$ replaced by $(u_k, \A _k)$,
where 
$\eta_k(x) = k^3 \eta (kx)$, $\eta \in\mathcal{S}$.
Let $v_k$ be the corresponding solution to the regularized equation
mentioned above.
Then $v_k$ satisfies the estimate \eqref{eq:H2estLS}.
Accordingly $\sup_k \| v_k; L^\infty(I;H^2) \| <\infty$
by virtue of Lemma \ref{nablaes}. 
Therefore there exists a subsequence of $\{ v_k\}$ 
that converges to some function $v \in L^\infty (I;H^2)$
in $\textrm{w}^\ast$-sense.
We can easily check that  $v$ satisfies \eqref{eq:H2estLS}.
%This subsequence is also written as $\{ v_k\}$ for simplicity.
The function $v$ belongs to $W^{1,1}(I;L^2)$
and satisfies \eqref{eq:LS} almost every $t$
since $v$ satisfies the integral version of 
\eqref{eq:LS}-\eqref{eq:LSinitial}, namely
\begin{equation}\label{eq:intLS}
v(t)=v_0 -i\int_{t_0}^t [\mathcal{H}(\A )+ \phi(u)]v(\tau) d\tau.
\end{equation}
For, each $v_k$ clearly satisfies \eqref{eq:intLS}
with $(u,\A )$ replaced by $(u_k,\A _k)$,
and $\{ v_k \}$ converges to $v$ in $\textrm{w}^\ast$-sense as $k\to \infty$ 
along some suitable subsequence.
Finally we prove the strong continuity of $v(t)$ in $H^2$.
%
%Since $v \in L^\infty (I;H^2) \cap W^{1,\infty}(I;L^2)$
%$v(t)$ is weakly continuous in $H^2$, and 
%
To this end we remark that 
$v \in C_{\rm w}(I;H^2) \cap C(I;H^s)$ with $s<2$ by \eqref{eq:intLS} and that
$\A  \in C(I;L^p)$ with $3\le p <6$ by \ref{A1}.
Hence we can show that 
$2i\A  \nabla v +|\A |^2 v$ 
is strongly continuous in $L^2$ and that 
$\mathcal{H}(\A (t))v(t)$ is weakly continuous in $L^2$.
We use the estimate \eqref{eq:H2estLS} with $I=[t_0,t]$
and the conservation of the $L^2$-norm to obtain
$\limsup_{t\to t_0} \| \mathcal{H}(\A (t))v(t) \|_2 
\le \| \mathcal{H}(\A (t_0))v(t_0) \|_2$. 
This inequality and the weak continuity conclude the strong continuity of 
$\mathcal{H}(\A (t)) v(t)$ in $L^2$, and hence 
$v(t)$ is strongly continuous in $H^2$.
The last part of the lemma is so easy that we omit the proof.
\ebox

\thmskip

By virtue of the lemma above, we can define the evolution operator 
for \eqref{eq:LS}.
Under the assumption \ref{A1},
we define a two-parameter family of operators 
$\{ U_{u,\A }(t,\tau)\}_{t,\tau \in I}$ by the relation
\begin{equation}\label{eq:LSevol}
U_{u,\A }(t,\tau) v(\tau) = v(t).
\end{equation}
Namely, we arbitrarily give the initial data at the time $\tau$, 
say $v(\tau)$, and solve \eqref{eq:LS} up to the time $t$;
then we define the image of $v(\tau)$ by $v(t)$.
In what follows we omit the lower indices $u,\A $ 
unless it causes any confusion.
Clearly this family of operators is well-defined,
and has the group property:
\begin{equation}
\begin{split}\label{eq:group}
U(t,\tau)U(\tau,\tau')&=U(t,\tau'), \\
U(t,t)&=\I
\end{split}
\end{equation}
for $t,\tau\in I$.
On account of Lemma \ref{lem:LSH2}
together with Lemma \ref{equinorm},
$U(t,\tau)$ are uniformly bounded operators on $H^2$ with the estimate
\begin{equation}
\label{eq:equiconti}
\begin{split}
K_2 &\equiv \sup_{t,\tau \in I} \| U(t,\tau) ;H^2 \to H^2 \| \\
&\lesssim \{ 1+\| \A ;L^\infty (I;H^1)\| \}^4 
\exp \bigl\{ C \| u;L^2(I;H^{3/4}) \|^2 
  +C\| \partial_t \A  ; L^1(I;L^3)\| \bigr\}.
\end{split}
\end{equation}
This family is strongly continuous in $H^2$.
Namely, for any $\psi \in H^2$, the function
\[ (t,\tau)\in I\times I \mapsto U(t,\tau)\psi \]
is strongly continuous in $H^2$.
Indeed, $U(t,\tau)\psi$ is strongly continuous in $t$. 
Combining this fact with \eqref{eq:group}-\eqref{eq:equiconti},
we obtain the strong continuity as a two-variable function.
\begin{lemma}\label{lem:LSHm}
Under the assumption {\rm \ref{A1},}
$\{ U(t, \tau)\}$ defined by \eqref{eq:LSevol}
can be uniquely extended to a strongly continuous 
two-parameter family of operators on $H^s${\rm ,} 
$0\le s\le 2${\rm ,} with the estimate 
\begin{equation}\label{eq:HmestLS}
K_s \equiv \sup_{t,\tau \in I} \| U(t,\tau); H^s \to H^s \| \le K_2^{s/2}.
\end{equation}
Especially{\rm ,} $\{U(t,\tau)\}$ is a unitary group on $L^2$ and 
\begin{equation}\label{eq:unitary}
U(t,\tau)^\ast = U(\tau,t).
\end{equation}
Moreover{\rm ,} for any $v_0 \in H^s${\rm ,} $U(t,t_0)v_0$ is 
a unique $H^s$-solution to \eqref{eq:LS}-\eqref{eq:LSinitial}{\rm .}
\end{lemma}
\Proof
$\{U(t,\tau)\}$ can be extended as a family of unitary operators in $L^2$ 
on account of the $L^2$-norm conservation law and the fact 
that each $U(t,\tau)$ is a bijection on $H^2$.
Therefore \eqref{eq:HmestLS} is proved by interpolation.
Therefore the first part of the lemma has been proved except
strong continuity of $U$; 
this is a consequence of the continuity of the $H^s$-solution,
which is proved below.
The relation \eqref{eq:unitary} follows 
from the unitarity and the group property.
The latter part is proved by approximation.
Let $\{ v_{0j}\}_{j=1}^\infty \subset H^2$ be a sequence 
converging to $v_0$ in $H^s$.
Then $v_j(t)=U(t,t_0) v_{0j}$ are $H^2$-solutions
and the sequence $\{ v_j \}\subset C(I;H^2)$ strongly converges to 
$v(t) \equiv U(t,t_0) v_0$ in $L^\infty (I;H^s)$ 
by virtue of the estimate \eqref{eq:HmestLS}.
Hence $v \in C(I;H^s)$.
Moreover, $v$ satisfies \eqref{eq:intLS}
since each $v_j$ satisfies this equation 
with $v_0$ replaced by $v_{0j}$,
and since $v_j \to v$ strongly in $C(I;H^s)$.
This fact implies that $v \in W^{1,1}(I;H^{s-2})$.
Therefore $v$ is an $H^s$-solution.
Finally we show uniqueness, 
which has yet to be proved in the case $s<1$.
Let $v(t)$ be an $H^s$-solution, and $\psi \in H^2$ 
be an arbitrary function.
Then
\begin{align*}
\frac{d}{dt}\langle v(t),U(t,\tau)\psi \rangle
=\langle -i(\mathcal{H}+\phi)v, U(t,\tau)\psi \rangle
  +\langle v(t),-i(\mathcal{H}+\phi)U(t,\tau)\psi \rangle
=0.
\end{align*}
Therefore
$\langle v_0,U(t_0,\tau)\psi \rangle
=\langle v(t_0),U(t_0,\tau)\psi \rangle
=\langle v(\tau),U(\tau,\tau)\psi \rangle
=\langle v(\tau),\psi \rangle$.
This means $v(\tau)=U(\tau,t_0)v_0$ for any $\tau\in I$.
%since $U(\tau,t_0)=U(t_0,\tau)^\ast$ 
%by the $L^2$-norm conservation law.
Therefore the uniqueness has been proved.
\ebox

\begin{corollary}\label{cor:LSHm}
Under the assumption {\rm \ref{A1},}
$\{ U(t,\tau) \}$ defined by \eqref{eq:LSevol} can be extended 
uniquely to a strong continuous family on $H^{-s}${\rm ,} 
$0<s\le 2${\rm ,} and 
\[ \| U(t,\tau); H^{-s}\to H^{-s}\|\le K_s. \] 
\end{corollary}
\Proof
The corollary immediately follows from Lemma \ref{lem:LSHm}
by duality.
\ebox

\thmskip

Next we consider the inhomogeneous problem.

\begin{lemma}\label{lem:inhomogeneous}
We assume {\rm \ref{A1}.}
Let $f\in L^1(I;H^{-2})${\rm,} and
$v\in C(I;L^2)\cap W^{1,1}(I;H^{-2})$ be an $L^2$-solution to 
\begin{equation}\label{eq:inhomogeneous}
i\partial_t v =\mathcal{H}(\A )v +\phi(u)v +f.
\end{equation}
Then for any $t_0\in I${\rm ,} 
\begin{equation}\label{eq:intinhom}
v(t)=U(t,t_0) v(t_0)-i\int_{t_0}^t U(t,\tau)f(\tau)d\tau.
\end{equation} 
Here $\{ U(t,\tau) \}$ is the evolution operator for \eqref{eq:LS}{\rm .}
\end{lemma}
\Proof
We take $\psi \in H^2$ arbitrarily.
Then $\langle v(\tau), U(\tau,t)\psi \rangle$ is absolutely continuous 
with respect to $\tau$ and for almost every $\tau\in I$
\[ \frac{d}{d\tau}\langle v(\tau), U(\tau,t)\psi \rangle
= \langle -if(\tau), U(\tau,t)\psi \rangle
= -i \langle U(t,\tau)f(\tau), \psi \rangle. \]
Integrating this formula with respect to $\tau$ on $[t_0,t]$, we obtain
\[ \langle v(t), \psi \rangle
=\langle U(t,t_0) v(t_0), \psi \rangle
-i \int_{t_0}^t \langle U(t,\tau)f(\tau) , \psi \rangle d\tau. \]
This means \eqref{eq:intinhom}.
\ebox

\thmskip
We proceed to the case $s>2$.
\begin{lemma}\label{lem:HsestLS}
Let $s>2${\rm ,} $\sigma \ge \max\{s-2,(2s-1)/4,1\}$ and $(s,\sigma) \neq (7/2,3/2)${\rm .}
Let $\A \in \bigcap_{j=0}^2 W^{j,\infty}(I;H^{\sigma-j})$ 
with $\Div \A =0$ and $\partial_t \A \in L^1(I;L^3)${\rm .}
Let $u \in \bigcap_{j=0}^2 W^{j,\infty}(I;H^{s-2j})${\rm .}
Let $v_0\in H^s$. 
Then the $H^2$-solution $v$ to \eqref{eq:LS}-\eqref{eq:LSinitial}
actually belongs to $C_{\rm w}(I;H^s)$ and satisfies
\begin{equation}\label{eq:HsestLS}
\begin{split} 
\| v ;L^\infty (I;H^s)\| &\lesssim
K_{s-4}\langle M^1_\sigma \vee R^1_{s-1} \rangle^\alpha \| v_0 ;H^s\| \\
&\qquad\times \exp\left(
CK_{s-4} \langle M^2_\sigma \vee R^2_s \rangle^4
\int_I \langle \| \partial_t \A \|_3 \rangle dt
\right).
\end{split}
\end{equation}
Here 
$M^k_\sigma = \max_{0\le j \le k} \| \partial_t^j \A; L^\infty(I;H^{\sigma -j}) \|${\rm ,}
$R^k_s = \max_{0\le j \le k} \| \partial_t^j u; L^\infty(I;H^{s-2j}) \|$
and $\alpha$ is some positive constant{\rm .}
Moreover if $v=u$ and $s\ge5/2${\rm ,} we have 
\begin{equation}\label{eq:HsestNLS}
\begin{split} 
\| u ;L^\infty (I;H^s)\| &\lesssim
K_{s-4}\langle M^1_\sigma \vee R^1_{s-1} \rangle^\alpha \| u_0 ;H^s\| \\
&\qquad\times \exp\left(
CK_{s-4} \langle M^2_\sigma \vee R^2_{s-1} \rangle^\alpha
\int_I \langle \| \partial_t \A \|_3 \rangle dt
\right).
\end{split}
\end{equation}

\end{lemma}

\begin{remark}
For $-2\le s\le 2$, $K_s$ is defined and estimated as in Lemma \ref{lem:LSHm} and Corollary \ref{cor:LSHm}.
Therefore once we have obtained \eqref{eq:HsestLS},
this estimate ensures that $\{U(t,\tau)\}$ is a family of operators on $H^s$
and gives an upper bound of $K_s$ for $2<s\le 6$.
Repeating this process, we can inductively estimate $K_s$ for all $s>2$.
\end{remark}

\Proof 
In the following proof, the exponent $\alpha$ may be different line to line;
precisely we have to replace $\alpha$ by the greatest one that has ever appeared,
but for simplicity we omit this process and use the same letter $\alpha$.
We estimate $\| \partial_t^2 u ;H^{s-4}\|$ instead of $\| u ;H^s \|$
since they are expected to be equivalent.
To this end, we differentiate \eqref{eq:LS} in $t$ twice.
By simple calculation, we obtain
\begin{align}
i\partial_t^2v &= (\calh+\phi)\partial_tv+(2i\partial_t\A(\nabla-i\A)+\partial_t\phi)v,
\label{v2} \\
i\partial_t^3v &= (\calh+\phi)\partial_t^2v+4i\partial_t\A(\nabla-i\A)\partial_t v \notag \\
   &\quad +2\partial_t\phi\partial_tv+2i\partial_t^2\A(\nabla-i\A)v+2(\partial_t\A)^2v+\partial_t^2\phi v \notag\\
&\equiv (\calh+\phi)\partial_t^2v+F_1+\cdots+F_5.
\label{v3}
\end{align}
We put $F\equiv\sum_{j=1}^5F_j$. By virtue of Lemma \ref{lem:inhomogeneous}, we convert \eqref{v3} to the integral form.
Precisely we need $s\ge 4$ to apply the lemma, but we have the expression below for $s>2$ by regularizing technique:
\begin{equation}\label{eq:v3int}
\partial_t^2v(t)=U(t,t_0)\partial_t^2v(t_0)-i\int_{t_0}^tU(t,\tau)F(\tau)d\tau.
\end{equation}
Therefore 
\begin{equation}\label{eq:s-4est}
\|\partial^2_tv(t);H^{s-4}\|\le K_{s-4}\{\|\partial_t^2v(t_0);H^{s-4}\|+\int_{t_0}^t\|F(\tau);H^{s-4}\|d\tau\}.
\end{equation}
We estimate the right-hand side.
First, we prove that the following equivalence holds for the solution $v$:
\begin{equation}\label{eq:equiv}
\|v;H^s\| +\langle M^1_\sigma \vee R^1_{s-1} \rangle^\alpha \| v\|_2
\simeq
\|\partial_t^2v;H^{s-4}\|+\langle M^1_\sigma \vee R^1_{s-1} \rangle^\alpha \| v\|_2.
\end{equation}
By the use of \eqref{equi2} twice and equations \eqref{eq:LS}, \eqref{v2}, we obtain for $s>2$
\begin{align*}
\|v;H^s\|&\lesssim \|\partial_t^2v;H^{s-4}\|+\|\partial_t\A(\nabla-i\A)v;H^{s-4}\| +\|\partial_t\phi v;H^{s-4}\| \\
&\quad +\langle M^0_\sigma \vee R^0_{s-1} \rangle^\alpha\|\partial_tv\|_2 
+\langle M^0_\sigma \vee R^0_{s-1} \rangle^\alpha\|v\|_2\\
&\equiv \|\partial_t^2v;H^{s-4}\|+F_6+\cdots+F_9.
\end{align*}
We begin with the estimate of 
$F_6$.
If $s\ge 4$, by the Leibniz formula and \eqref{nablaes1} we have 
\begin{equation}
\begin{split}
&\|\partial_t\A(\nabla-i\A)v;H^{s-4}\| \\
&\quad\lesssim 
\|\partial_t\A;H^{s-4,6}\|\|(\nabla-i\A)v\|_3+\|\partial_t\A\|_6\|(\nabla-i\A)v;H^{s-4,3}\| \\
&\quad\lesssim \|\partial_t\A;H^{s-3}\|\|(\nabla-i\A)v;H^{s-3}\|  \\
&\quad\lesssim \|\partial_t\A;H^{s-3}\|\|v;H^{s-2}\|\langle\|\A;H^\sigma\|\rangle.
\end{split}\label{dAnv}
\end{equation}
If $2<s<4$, we choose $p$ so that 
$(1/2-(s-2)/3)_+\le 1/p \le 1/2-(1/2-(4-s)/3)_+$.
Then by the Sobolev inequality  we have the continuous embeddings 
$H^{s-2}\embedding L^p$ 
and
$L^q\embedding H^{s-4}$, where $1/q=1/2+1/p$.
Using these embeddings together with \eqref{nablaes1} we obtain
\begin{align*}
&\|\partial_t\A(\nabla-i\A)v;H^{s-4}\| \lesssim \|\partial_t\A(\nabla-i\A)v \|_q 
\lesssim \|\partial_t\A\|_2\|(\nabla-i\A)v\|_p  \\
&\quad\lesssim \|\partial_t\A\|_2\|(\nabla-i\A)v;H^{s-2}\| 
\lesssim \|\partial_t\A\|_2\|v;H^{s-1}\|\langle\|\A;H^\sigma\|\rangle.
\end{align*}
Therefore we get
\[ F_6 \lesssim \langle M^1_\sigma \rangle^2 \| v; H^{s-1}\|. \]
Next we derive the estimate for $F_7$. We have
\[ F_7 \lesssim  \langle R^1_{s-1} \rangle^2 \|v;H^{s-2}\|. \]
If $s\ge 3$, 
this is proved by Lemma \ref{phies} as 
\[ F_7 \lesssim \|\partial_t \phi u;H^{s-3}\|\lesssim \|\partial_t u;H^{s-3}\| \|u;H^{s-1}\|\|v;H^{s-2}\|. \]
If $2<s<3$, we use the duality estimate as follows,
from which we obtain the desired result:
\begin{align*}
|\langle\omega^{-2}(\partial_tu\bar{u})v,\psi\rangle| 
&= |\langle\partial_t u,\omega^{-2}(\bar{v}\psi)u\rangle| 
\le \|\partial_tu;H^{s-3}\|\|\omega^{-2}(\bar{v}\psi)u ;H^{3-s}\| \\
&\lesssim \|\partial_tu;H^{s-3}\|\|v;H^{s-2}\|\|\psi;H^{4-s}\|\|u;H^{s-1}\|.
\end{align*}
The estimate for $F_8$ is easy. Indeed we obtain by \eqref{nablaes1} 
\[ F_8 \lesssim \|v;H^2\|\langle\|\A;H^\sigma\|\vee\|u\|_2\rangle^\alpha
\lesssim \langle M^1_\sigma \vee R^0_{s-1} \rangle^\alpha \|v;H^2\|.\]
Using these estimates with \eqref{interpolation}, we obtain for any $\varepsilon>0$
\begin{equation}
\sum_{j=6}^9 F_j \le 
C(\varepsilon)\|v\|_2\langle M^1_\sigma \vee R^1_{s-1} \rangle^\alpha
+\varepsilon\|v;H^s\|,
\label{remainder}
\end{equation}
where $C(\varepsilon)$ is a positive constant.
Therefore we have proved
\[ 
\|v;H^s\| \lesssim
\|\partial_t^2v;H^{s-4}\|+\langle M^1_\sigma \vee R^1_{s-1} \rangle^\alpha \| v\|_2.
\]
The opposite inequality in \eqref{eq:equiv} is similarly proved.
Applying this inequality to \eqref{eq:s-4est} together with the $L^2$-norm conservation law, we obtain 
the following intermediate estimate:
\[ \|v(t);H^s\| \lesssim K_{s-4}\{ 
			\langle M^1_\sigma \vee R^1_{s-1} \rangle^\alpha \| v_0; H^s \|
			+\int_{t_0}^t \| F(\tau);H^{s-4}\| d\tau \}. \]

The next step is the estimate of $F=\sum_{j=1}^5 F_j$ in $H^{s-4}$.
We begin with the estimate of $F_1$.
If $s\ge 4$, by the estimate \eqref{dAnv} with $v$ replaced by $\partial_t v$, we have 
\[ \| F_1 ;H^{s-4} \| \lesssim \|\partial_t\A;H^{s-3}\| \langle\|\A;H^\sigma\|\rangle \|\partial_t v;H^{s-2}\|.
%\lesssim \langle M^1_\sigma \rangle^2 \|\partial_t v;H^{s-2}\|.
\]
If $3\le s<4$, we choose $1/p=1/2+(4-s)/3$, $1/q=1/2-(s-3)/3$.
Then by the Sobolev inequality and \eqref{nablaes1}
\begin{align*}
\|F_1;H^{s-4}\| &\lesssim \|\partial_t\A (\nabla-i\A)\partial_tv\|_p 
\lesssim \|\partial_t\A\|_3\|(\nabla-i\A)\partial_tv \|_q \\
&\lesssim \|\partial_t\A\|_3\|(\nabla-i\A)\partial_tv;H^{s-3}\| 
\lesssim \|\partial_t\A\|_3\langle \|\A;H^1\|\rangle\|\partial_tv;H^{s-2}\|.
\end{align*}
For $2<s<3$, we use $\Div \A=0$ and \eqref{nablaes1} to have
\begin{align*}
\|F_1;H^{s-4}\| &= \|(\nabla-i\A)(\partial_t\A\partial_tv);H^{s-4}\| 
\lesssim \|\partial_t\A\partial_tv;H^{s-3}\|\langle\|\A;H^1\|\rangle \\
&\lesssim \|\partial_t\A\|_3\langle\|\A;H^1\|\rangle\|\partial_tv;H^{s-2}\|.
\end{align*}
Here we have used the Sobolev inequality twice as in the previous case.
In any cases we have by \eqref{Hphi}
\[ \|F_1;H^{s-4}\|\lesssim \langle M^0_\sigma \vee R^0_{s-1}\rangle^4
\langle\|\partial_t\A\|_3\rangle \| v;H^s\|. \]
By \eqref{Hphi} and the estimate for $F_7$ with $v$ replaced by $\partial_t v$, we have 
\[
\|F_2;H^{s-4}\| \lesssim \langle M^0_{\sigma}\vee R^1_{s-1}\rangle^4 \|v;H^s\|.
\]
The estimate of $F_3$. 
If $s\ge 4$, by the Leibniz formula and \eqref{nablaes1} we have 
\begin{align*}
\|F_3;H^{s-4}\| &\lesssim \|\partial_t^2\A;H^{s-4}\|\|(\nabla-i\A)v\|_\infty
+\|\partial_t^2\A;1/2-\varepsilon\|\|(\nabla-i\A)v;H^{s-4,\varepsilon^{-1}}\| \\
&\lesssim \|\partial_t^2\A;H^{s-4}\|\|(\nabla-i\A)v;H^{s-2}\| 
\lesssim \|\partial_t^2\A;H^{s-4}\|\langle\| \A ;H^\sigma \|\rangle \|v;H^{s-1}\|,
\end{align*}
where $\varepsilon$ is a sufficiently small number,
and the second term in the right-hand side of the first inequality is removed if $s=4$.
If $2<s<4$, we use duality.
By \eqref{nablaes2}, we have 
\begin{align*}
|\langle F_3,\psi\rangle| &= |\langle \partial_t^2\A,\psi(\nabla+i\A)\bar{v}\rangle| \\
&\le \|\partial_t^2\A;H^{\sigma_0-2}\|\|\psi;H^{4-s}\|\|v;H^{\sigma_0+1}\|
\langle\|\A;H^{\sigma_0}\|\rangle,
\end{align*}
with $1\vee(2s-1)/4\vee(s-2)\le \sigma_0\le 2$, $(s,\sigma_0)\neq(7/2,3/2), (5/2,1)$.
Taking $\sigma_0$ so that $\sigma_0\le \sigma\wedge(s-1)$, we obtain 
\[
\|F_3;H^{s-4}\| \lesssim \| \partial_t^2\A;H^{\sigma-2}\| \langle\| \A ;H^\sigma \|\rangle\|v;H^s\|.
\]
Therefore we have 
\[ F_3 \lesssim \langle M^2_\sigma \rangle^2 \|v;H^s\|. \]
The estimate of $F_4$. If $s\ge 4$, we have by the Leibniz rule and the Sobolev inequality
\begin{align*}
\|F_4;H^{s-4}\| &\lesssim \|\partial_t\A;H^{s-4,6}\|\|\partial_t\A\|_3 \|v\|_\infty
+\|\partial_t\A\|_6^2\|v;H^{s-4,6}\| \\
&\lesssim \|\partial_t\A; H^{\sigma -1}\|^2 \|v;H^s\|.
\end{align*}
If $s<4$, we have by the Sobolev inequality
\begin{align*}
\|F_4;H^{s-4}\|&\lesssim
\begin{cases}
\|\partial_t\A\|_3 \|\partial_t\A\|_2 \|v \|_\infty & \text{if } s\le 3, \\
\|\partial_t\A\|_3 \|\partial_t\A;1/2-(s-3)/3\|\|v \|_\infty & \text{if } s>3
\end{cases}\\
&\lesssim \|\partial_t\A;L^3\| \| \partial_t \A ;H^{\sigma-1}\| \|v;H^s\|.
\end{align*}
Therefore we have
\[
\|F_4;H^{s-4}\| \lesssim \langle \|\partial_t\A;L^3\|\rangle \langle M^0_\sigma\rangle^2 \|v;H^s\|.
\]
The estimate for $F_5$. 
If $s\ge 4$, by Lemma \ref{phies} we have 
\begin{align*}
& \|\omega^{-2}(\partial_t^2u\bar{u})v;H^{s-4}\|\lesssim \|\partial_t^2u;H^{s-4}\|\|u;H^{s-2}\|\|v;H^{s-4}\|, \\
& \|\omega^{-2}(\partial_tu\partial_t\bar{u})v;H^{s-4}\|\lesssim \|\partial_tu;H^{s-3}\|^2\|v;H^{s-2}\|,
\end{align*}
which lead 
\[
\|F_5;H^{s-4}\| \lesssim (R^2_s)^2\|v;H^{s-2}\|.
\]
If $2<s<4$, it is sufficient to show 
\begin{align*}
|\langle\omega^{-2}(\bar{u}\partial_t^2u)v,\psi\rangle| &=
|\langle \partial_t^2u,\omega^{-2}(\bar{v}\psi)u\rangle| \\
&\lesssim \|\partial_t^2u;H^{s-4}\|\|v;H^{s-1}\|\|\psi;H^{4-s}\|\|u;H^{s}\|, \\
\|\omega^{-2}(\partial_tu\partial_t\bar{u})v;H^{s-4}\| 
&\lesssim \|\partial_tu;H^{s-2}\|^2\|v;H^{s-1}\|,
\end{align*}
where we have used Lemma \ref{phies}, and $L^2\embedding H^{s-4}$ at the second inequality.
Therefore we obtain 
\[
\|F_5;H^{s-4}\| \lesssim \langle R^2_s\rangle^2\|v;H^{s-1}\|.
\]

Collecting all the estimates, we obtain 
\begin{align*}
\|v(t);H^s\| &\lesssim K_{s-4}\{ \|v_0;H^s\|\langle M^1_{\sigma}\vee R^1_{s-1}\rangle^\alpha \\
&\quad + \int_{t_0}^t \|v(\tau) ;H^s\| \langle \|\partial_t\A\|_3 \rangle
\langle M^2_{\sigma} \vee R^2_s\rangle^4 d\tau\},
\end{align*}
which leads \eqref{eq:HsestLS} by the Gronwall inequality.

To prove \eqref{eq:HsestNLS}, we only have to modify the estimate of $F_5$ for $5/2\le s<4$.
Indeed, if $s\ge 5/2$ we have
\begin{align*}
\|\omega^{-2}(\bar{u}\partial_t^2u)u ;H^{s-4}\| 
&\lesssim \|\partial_t^2u;H^{s-4}\| \|u;H^{s-1}\|^2 \\
&\lesssim \langle \| A;H^\sigma \|\vee \| u ;H^{s-1}\|\rangle^\alpha \|u;H^s\|
\end{align*}
and
\begin{align*}
\|\omega^{-2}(\partial_tu\partial_t\bar{u})u;H^{s-4}\| 
&\lesssim \|\partial_tu;H^{s-2}\| \|\partial_tu;H^{s-3}\| \|u;H^{s-1}\| \\
&\lesssim \langle \| A;H^\sigma \|\vee \| u ;H^{s-1}\|\rangle^\alpha \|u;H^s \|
\end{align*}
by virtue of the duality argument.
\ebox

\begin{corollary}\label{cor:Dtv}
Let $s, \sigma, u,\A$ and $v_0$ satisfy the assumption of Lemma {\rm \ref{lem:HsestLS}.}
Then for the solution $v$ to \eqref{eq:LS}-\eqref{eq:LSinitial}
we have the following.

{\rm (1)} The estimate
\[ \max_{j=1,2}\| \partial_t^j v ;L^\infty(I;H^{s-2j})\|
\lesssim \langle M^1_\sigma \vee R^1_{s-1}\rangle^\alpha \| v;L^\infty(I;H^s)\| \]
holds{\rm .} 
Here $M^1_\sigma, R^1_{s-1}$ are defined as in Lemma {\rm \ref{lem:HsestLS},}
and $\alpha$ is some positive number{\rm .}

{\rm (2)}
If $\A \in C(I;H^\sigma)${\rm ,} then
$v\in \bigcap_{j=0}^2 C(I;H^{s-2j})${\rm .}
Especially $\{ U(t,\tau)\}$ is a strongly-continuous family on $H^s${\rm .}
\end{corollary}
\Proof
(1) is the consequence of \eqref{Hphi} and \eqref{eq:equiv}.
We prove (2) for $2<s\le 6$. 
$F=\sum_{j=1}^5 F_j$ in \eqref{v3} belongs to $L^1(I;H^{s-4})$.
Therefore by virtue of Lemma \ref{cor:LSHm} together with Lebesgue's convergence theorem,
the right-hand side of \eqref{eq:v3int}, and hence $\partial_t^2 v$ belong to $C(I;H^{s-4})$.
To prove $\partial_t v\in C(I;H^{s-2})$, it suffices to show $\Delta \partial_t v \in C(I;H^{s-4})$
since we have $\partial_t v \in C(I;L^2)$.
If we recall the remark for Lemma \ref{equinorm} and use the estimates for $F_6, F_7$ 
in the proof of Lemma \ref{lem:HsestLS}, 
we can show that all the terms except $\Delta \partial_t v$ in \eqref{v2} belong to $C(I;H^{s-4})$. 
Therefore $\Delta \partial_t v \in C(I;H^{s-4})$.
Analogously we can prove $v\in C(I;H^s)$ by  \eqref{eq:LS}.
General case is proved by induction.
\ebox

\section{Linearized equation for {\bf MS-C}}

To solve \textbf{MS-C}, we consider the linearized equation below:
\begin{alignat}{2}
& i\partial_t v=(\calh(\A)+\phi(u))v, & \quad & v(0)=u_0,  \label{v}\\
& (\partial_t^2-\Delta+1)\B=P\J(u,\A)+\A, & \quad & B(0)=\A_0,\partial_t\B(0)=\A_1.\label{B}
\end{alignat}
We always assume $\Div \A=0$ and $(u_0,\A_0,\A_1)\in X^{s,\sigma}$.
In later sections we often use the equations with $(u,\A, v,\B)$ replaced by $(u',\A',v',\B')$ 
and $(u_0,\A_0,\A_1)$ by $(u_0',\A_0',\A_1')$.
We refer to such equations as $\eqref{v}', \eqref{B}'$, and
we often abbreviate $\calh(\A'), \phi(u'), \J(u',\A')$ to $\calh', \phi', \J'$.
If we define the map 
\[ \Phi: (u,\A)\mapsto (v,\B),\]
 the fixed points of $\Phi$ solve \textbf{MS-C}.
To prove the unique existence of the solution, 
we show that $\Phi$ is a contraction map in some appropriate function space.
In this section we prove a priori estimates for \eqref{v}-\eqref{B} which yields 
that $\Phi$ is a map from some function space to itself.
We treat \eqref{v} by the linear estimates discussed in Section 3,
and \eqref{B} mainly by the Strichartz estimate for the Klein-Gordon equation stated in Lemma \ref{wavest}.
We also need Lemma \ref{Pestimate} to treat the nonlinear term in \eqref{B}.

The Klein-Gordon equation
\begin{gather*}
 (\partial_t^2 -\Delta +1)A =f, \\
A(0)=A_0, \quad \partial_t A(0) = A_1
\end{gather*}
is solved as
\begin{equation}\label{eq:KG}
A(t) =\cos t\Omega A_0+\frac{\sin t\Omega}{\Omega}A_1 +\int_0^t \frac{\sin (t-\tau)\Omega}{\Omega}f(\tau)d\tau.
\end{equation}
We call a pair $(q,r)$ admissible if $2\le r<\infty$, $1/r+1/q=1/2$.
We put $\beta(r)\equiv 1-2/r=2/q$. 
With this notation, we have the following.

\begin{lemma} 
\label{wavest}
Let $(q_j,r_j)${\rm ,} $j=0,1${\rm ,} be any admissible pairs{\rm ,}
$I\subset\R$ be an interval containing $0${\rm .}
Let $(A_0,A_1) \in H^\sigma \oplus H^{\sigma-1}$ with $\sigma \in \R${\rm ,}
and $f \in L^{q'_1}(I;H^{\sigma -1+\beta(r_1),r'_1})${\rm .}
Then $A$ in \eqref{eq:KG} belongs to $C(I;H^\sigma) \cap C^1(I;H^{\sigma-1})$
and satisfies 
\begin{align*} 
&\| A ; L^{q_0}(I;H^{\sigma-\beta(r_0),r_0}) \| 
+ \| \partial_t A ; L^{q_0}(I;H^{\sigma-\beta(r_0)-1,r_0}) \| \\
&\quad\lesssim \| (A_0,A_1) ; H^\sigma \oplus H^{\sigma-1} \|
+\| f; L^{q'_1}(I;H^{\sigma -1+\beta(r_1),r'_1})\|. 
\end{align*}
\end{lemma}
\Proof
See for example \cite{B84,GV85AIHP,GV95,S77}.

%%%%%%%%%%%%%%%%%%%%%%%% Lemma {Pestimate}
\begin{lemma}
\label{Pestimate}
{\rm(1)} Let $s,\sigma,p$ satisfy $s>1, \sigma\ge0,\max\{3/(s-1),2\} \le p<\infty${\rm ,} 
$(s,p)\neq(5/2,2)${\rm .} 
Then
\[ \|P(u\nabla v); H^{\sigma,p'}\|\lesssim \|u;H^\sigma\|\|v;H^s\|+\|u;H^s\|\|v;H^\sigma\|. \]

{\rm(2)}
Moreover if $1\le \sigma\le s$ or $\sigma\le s-1${\rm ,} then 
\[
\|P(u\nabla v); H^{\sigma,p'}\| \lesssim \|u;H^\sigma\|\|v;H^s\|.
\]
\end{lemma}

\Proof 
We rewrite and estimate the left-hand side as 
\[
\|P(u\nabla v);H^{\sigma,p'}\|=\|P([\Omega^\sigma,u]\nabla v-\nabla u \Omega^\sigma v)\|_{p'} 
\lesssim \| [\Omega^\sigma,u]\nabla v-\nabla u \Omega^\sigma v \|_{p'},
\]
where we have used the property 
$P(u\nabla w)=-P(\nabla uw)$ and the fact that $P$ is a bounded operator on $L^{p'}$.
By the Kato-Ponce commutator estimate (see Appendix in \cite{KP88}, Lemma 2.10 in \cite{KPV91}), we have
\[
\|[\Omega^\sigma,u]\nabla v-\nabla u \Omega^\sigma v\|_{p'} 
\lesssim \| \Omega^\sigma u\|_2\|\nabla v\|_{2p/(p-2)}+\|\nabla u\|_{p_1}\| \Omega^\sigma v\|_{p_2}
\]
with $1/p'=1/p_1+1/p_2$, $p'\le p_2 <\infty$.
By putting $p_2=2$ and using the embedding $H^{s-1}\embedding L^{2p/(p-2)}$, we obtain the required 
estimate (1). 
The proof for (2) when $\sigma\le s-1$ follows from the direct application of 
the Leibniz formula to $\|u\nabla v;H^{\sigma,p'}\|$ and the Sobolev inequality.
For the proof of (2) when $1\le \sigma\le s$, we may assume $s-1<\sigma\le s$.
Putting $1/p_2=1/2-(s-\sigma)/3$, we obtain
$\|\nabla u\|_{p_1}\| \Omega^\sigma v\|_{p_2} \lesssim \|u;H^\sigma\|\|v;H^s\|$
again by the Sobolev inequality.
\ebox

%%%%%%%%%%%%%%%%%%% Lemma {W}
\begin{lemma}
\label{W}
Let $s,\sigma$ satisfy $s>1$ and $1\le \sigma\le \min\{(5s-2)/3,s+1\}$ with $(s,\sigma)\neq(5/2,7/2)${\rm .}
Let $\A,\B$ satisfy \eqref{B}{\rm .}
Then for $I=[0,T]$ with $0<T\le 1$ and
for any admissible pair $(q,r)$ the following estimate holds{\rm .}
\begin{equation}\label{eq:W1}
\begin{split}
& \max_{j=0,1}\|\partial_t^j\B;L^\infty(I;H^{\sigma-j})\cap L^q(I;H^{\sigma-\beta(r)-j,r})\| \\
&\quad\lesssim \|(\A_0,\A_1);H^\sigma\oplus H^{\sigma-1}\|
+T^{1/2}\langle\|u;L^\infty (I;H^s)\|\vee \|\A;L^\infty (I;H^{\sigma})\|\rangle^3.
\end{split}
\end{equation}
$\B \in \bigcap_{j=0}^1 C^j(I;H^{\sigma-j})$ if the right-hand side is finite{\rm .}
Moreover we have
\begin{equation}
\begin{split}\label{eq:W2}
&\|\partial_t^2\B;L^\infty(I;H^{\sigma-2})\| \\
&\quad\lesssim \|(\A_0,\A_1);H^\sigma\oplus H^{\sigma-1}\|
+\langle\|u;L^\infty (I;H^s)\|\vee \|\A;L^\infty (I;H^{\sigma})\|\rangle^3.
\end{split}
\end{equation}

\end{lemma}

\Proof 
By Lemma \ref{wavest}, the left-hand side of \eqref{eq:W1} is estimated by 
\begin{align*}
& \|(\A_0,\A_1);H^\sigma\oplus H^{\sigma-1}\|+\|\A;L^1 H^{\sigma-1}\| \\
&\quad+\|P(\bar{u}\nabla u);L^{q'_1}H^{\sigma+\beta(r_1)-1, r_1'}\| +\|\A|u|^2;L^1 H^{\sigma-1}\|,
\end{align*}
where $(q_1,r_1)$ is an admissible pair.
Under the assumption of the lemma,
there exists an exponent $r_1$ such that $\max\{2,3/(s-1)\} \le r_1<\infty$ and
$0\le \sigma+\beta(r_1)-1\le s$.
By (1) of Lemma \ref{Pestimate}, we have 
\[
\|P(\bar{u}\nabla u);H^{\sigma+\beta(r_1)-1,r_1'}\| 
\lesssim \|u;H^{\sigma+\beta(r_1)-1}\|\|u;H^s\| \lesssim \|u;H^s\|^2
\]
for this $r_1$.
On the other hand we have 
\begin{align*}
\|\A |u|^2;H^{\sigma-1}\| &\lesssim \|\A;H^{\sigma-1,6}\|\|u\|_6^2+
\|\A\|_\nu \|u;H^{\sigma-1,p_1}\|\|u\|_{p_2}\\
&\lesssim \|\A;H^\sigma\| \|u;H^s\|^2,
\end{align*}
where $1/\nu=1/2-1/p_1-1/p_2$. We choose $p_1=p_2=6$ if $\sigma\le s$, 
$1/p_1=1/2-(s+1-\sigma)/3$, $1/p_2=(1/2-s/3)_+$ if $s<\sigma\le s+1$ 
so that $H^s\embedding H^{\sigma-1,p_1}, L^{p_2}$.
With such a choice $H^\sigma \embedding L^\nu$ under the assumption of the lemma.
With these estimates and the H\"older inequality for the time variable, we obtain \eqref{eq:W1}.
Finally we prove \eqref{eq:W2}. By \eqref{B} we have
\[ \| \partial_t^2 \B ;L^\infty H^{\sigma-2}\|\le \| \B;L^\infty H^\sigma \|+\| \A+P\J ;L^\infty H^{\sigma-2}\|. \]
We have 
$\| P\J ;L^\infty H^{\sigma-2}\| \lesssim \langle\|u;L^\infty (I;H^s)\|\vee \|\A;L^\infty (I;H^{\sigma})\|\rangle^3$
similarly as above, since $H^{\sigma+\beta(r_1)-1,r_1'}\embedding H^{\sigma-2}$.
Therefore we obtain the required result.
\ebox

\thmskip

Now we define the function spaces where we consider the map $\Phi$.
We put for $s \le 2$
\begin{equation}\label{eq:space2}
\begin{split}
Z_{s,\sigma}&=\biggl\{ (u,\A) \in L^\infty(I;H^s\oplus H^\sigma); 
\| u;L^\infty(I;H^s)\|\le l_S, \\
&\qquad \Div\A=0, \A \in W^{1,6}(I;L^3), 
\| \A;L^\infty(I;H^\sigma)\|\vee\|\partial_t \A;L^6(I;L^3) \|\le l_M
\biggr\}.
\end{split}
\end{equation}
Here $I=[0,T]$. For $s>2$, we put 
\begin{equation}\label{eq:spaces}
Z_{s,\sigma}=Z_{s_\ast,\sigma_\ast} \cap \Dot{Z}_{s,\sigma}
\end{equation}
with
\begin{align*}
\Dot{Z}_{s,\sigma} &= \biggl\{ (u,\A) \in \bigcap_{j=0}^2 W^{j,\infty}(I;H^{s-2j} \oplus H^{\sigma-j});
\max_{0\le j\le 2}\| \partial_t^j u;L^\infty(I;H^{s-2j})\|\le L^s_S, \\
&\qquad\max_{j=0,1}\| \partial_t^j \A;L^\infty(I;H^{\sigma-j})\|\le L^\sigma_M,
\| \partial_t^2 \A;L^\infty(I;H^{\sigma-2})\|\le \tilde{L}^{\sigma}_M
\biggr\}.
\end{align*}
Here $s_\ast = (s-1)\vee 2$ and  $\sigma_\ast\le \sigma$ is a number such that 
$(s_\ast,\sigma_\ast)$ satisfies the assumption of Proposition \ref{prop:itself} below
with $(s,\sigma)$ replaced by $(s_\ast,\sigma_\ast)$.

\begin{proposition}\label{prop:itself}
Let $s\ge6/5${\rm ,} $\max\{4/3,s-2,(2s-1)/4\} \le \sigma\le \min\{s+1, (5s-2)/3\}$
and $(s,\sigma)\neq (5/2,7/2), (7/2,3/2)${\rm .}
Let the map $\Phi$ be defined by \eqref{v}{\rm ,} \eqref{B}{\rm .}

{\rm(1)}
If $s\le 2${\rm ,} there exist $l_S,l_M,T$ so that $\Phi$ is a map from $Z_{s,\sigma}$ to itself{\rm .}

{\rm(2)}
Let $s> 2$ and let $\Phi$ map $Z_{s_\ast,\sigma_\ast}$ to itself{\rm .}
Then there exist $L^s_S,L^\sigma_M,\tilde{L}^\sigma_M,T$ so that $\Phi$ is a map from $Z_{s,\sigma}$ to itself{\rm .}
\end{proposition}
\begin{remark}
If $\Phi(Z_{s,\sigma})\subset Z_{s,\sigma}$, then
the same inclusion holds even if we replace $T$ by a smaller one.
\end{remark}
\Proof
(1)
Let $(v,\B)=\Phi (u,\A)$.
By Lemmas \ref{lem:LSHm} and \ref{W} we have
\[ \| v;L^\infty (I;H^s) \|\le C \langle l_M\rangle^{2s} \exp (CTl_S^2+CT^{5/6}l_M) \| u_0 ;H^s \| \]
and
\[ \| \B;L^\infty (I;H^\sigma) \| \vee \| \partial_t \B;L^6 (I;L^3) \|
\le C\| (\A_0, \A_1);H^\sigma \oplus H^{\sigma-1} \| +CT^{1/2} \langle l_S \vee l_M \rangle^3. \]
We choose $l_S$, $l_M$ and $T$ as follows.
First, we choose $l_M$ so that 
$C\| (\A_0, \A_1);H^\sigma \oplus H^{\sigma-1} \| \le l_M/2$.
Next we choose $l_S$ so that 
$C \langle l_M\rangle^{2s} \| u_0 ;H^s \| \le l_S/2$.
Finally we choose $T$ so that
$\exp (CTl_S^2+CT^{5/6}l_M)\le 2$ and that $CT^{1/2} \langle l_S \vee l_M \rangle^3 \le l_M/2$.
Then we have 
$\| v;L^\infty (I;H^s) \|\le l_S$ and 
$\| \B;L^\infty (I;H^\sigma) \| \vee \| \partial_t \B;L^6 (I;L^3) \| \le l_M$.
Therefore $\Phi(Z_{s,\sigma})\subset Z_{s,\sigma}$.

(2)
By Lemma \ref{lem:HsestLS} together with Corollary \ref{cor:Dtv} and Lemma \ref{W} we have
\begin{align*}
\max_{0\le j \le 2} \| \partial_t^j v; L^\infty (I;H^{s-2j}) \|
&\le CK_{s-1}\langle L^\sigma_M\vee L_S^{s-1}\rangle^\alpha \\
&\qquad\times\exp (CK_{s-1}\langle L^\sigma_M \vee \tilde{L}^\sigma_M 
\vee L_S^s\rangle^4 T^{5/6} \langle l_M\rangle) \| u_0 ;H^s \|, 
\end{align*}
\[ \max_{j=0,1} \|\partial_t^j\B;L^\infty (I;H^{\sigma-j})  \|
\le C\| (\A_0, \A_1);H^\sigma \oplus H^{\sigma-1} \| +CT^{1/2} \langle L^s_S \vee L^\sigma_M \rangle^3 
\]
and 
\[ \|\partial_t^2\B;L^\infty (I;H^{\sigma-2})  \|
\le C\| (\A_0, \A_1);H^\sigma \oplus H^{\sigma-1} \| + \langle L^s_S \vee L^\sigma_M \rangle^3.
\] 
Since $K_{s-1}$, $l_M$ and  $L^{s-1}_S$ are increasing functions of $T$, 
we can choose suitable $L^\sigma_M,L^s_S,\tilde{L}^\sigma_M$ and $T$
in this order so that $\Phi(Z_{s,\sigma})\subset Z_{s,\sigma}$.
\ebox

\thmskip

We also need the following space: 
\begin{equation}\label{eq:space5/3}
\begin{split}
\tilde{Z}_{s,\sigma}&=\{ (u,\A) \in Z_{s,\sigma}\cap W^{1,\infty}(I;H^{s-2}\oplus H^{\sigma-1});
\partial_t \A \in C(I;L^2)\\
&\qquad \| \partial_t u;L^\infty (I;H^{s-2})\|\le l_S, \| \partial_t A;L^\infty (I;H^{\sigma-1})\|\le l_M \}.
\end{split}
\end{equation}
This space is mainly used to discuss the uniqueness for $s<7/4$.

\begin{proposition}\label{prop:itself5/3}
Let $6/5\le s \le 2$ and $4/3 \le \sigma\le (5s-2)/3${\rm .}
Then there exist $l_S,l_M$ and $T$ so that $\Phi$ is a map from $\tilde{Z}_{s,\sigma}$ to itself{\rm .}
\end{proposition}
\noindent
The proof is similar to that of the previous proposition, and left to the reader.

%%%%%%%%%%%%%%%%%%%%% The contraction argument for \do{s\ge7/4}
\section{The contraction argument for $s\ge7/4$}
We give the proof of Theorem \ref{thmsc} in the case of $s\ge7/4$.
We consider the metric
\begin{equation}
d(u,\A,u',\A')\equiv \|u-u';L^\infty (I;L^2)\|\vee 
\|\A-\A';L^\infty (I;H^{1/2})\cap L^4(I;L^4)\|. \label{eq:d} 
\end{equation}
We prepare the following proposition on this metric.
%%%%%%%%%%%
\begin{proposition}
\label{metric}
Let $I=[0,T]$ with $0<T\le 1${\rm .}
Let $(u,\A), (u',\A')\in Z_{7/4,4/3}${\rm .}
Let $(v,\B)$ and $(v',\B')$ be the solutions to 
\eqref{v}-\eqref{B} and $\eqref{v}'$-$\eqref{B}'$ respectively{\rm .}
Then the estimate 
\begin{align*}
d(v,\B,v',\B')&\lesssim \|(u_0-u_0',\A_0-\A_0',\A_1-\A_1');X^{0,1/2}\| \\
&\quad+T^{1/2}\langle \|v';L^\infty (I;H^{7/4})\| \vee l_S\vee l_M \rangle^2d(u,\A,u',\A')
\end{align*}
holds{\rm .}
%\[
%L\equiv \max_{w=u,u',v'}\|w;L^\infty (I;H^{7/4})\|\vee \max_{\Q=\A,\A'}\|\Q;L^\infty (I;H^{3/4})\|.
%\]
\end{proposition}

\Proof
We write the difference of the equations for $v$ and $v'$ as
\begin{equation}
i\partial_t(v-v')=(\calh+\phi)(v-v')+(\calh+\phi-\calh'-\phi')v'.
\label{1difv}
\end{equation}
We regard $(\calh+\phi-\calh'-\phi')v'$ as the inhomogeneous term and convert this equation
to the integral form by virtue of Lemma \ref{lem:inhomogeneous}.
%\[
%(v-v')(t)=U(t,0)(u_0-u_0')-i\int_0^tU(t,\tau)(\calh+\phi-\calh'-\phi')v'(\tau)d\tau.
%\]
%by virtue of Lemma \ref{lem:inhomogeneous}.
%Here $U(t,\tau)=U_{u,\A}(t,\tau)$. 
Since $U_{u,\A}(t,\tau)$ is unitary on $L^2$, we have 
\begin{equation}
\|(v-v')(t)\|_2\le \|u_0-u_0'\|_2+\|(\calh+\phi-\calh'-\phi')v';L^1L^2\|.
\label{difv}
\end{equation}
For the second term of the right-hand side we use the identity
\[ (\calh+\phi-\calh'-\phi')v'=2i(\A-\A')(\nabla-i(\A+\A')/2) v'+\omega^{-2}\textrm{Re}(u-u')(\overline{u-u'}) v'\]
%We write 
%\[
%\begin{split}
%\calh+\phi-\calh'-\phi'&=2i(\A-\A')\nabla+(\A+\A')(\A-\A')\\
%&\quad+\omega^{-2}((u-u')\bar{u}) +\omega^{-2}(u'(\overline{u-u'}))
%\end{split}
%\]
and estimate it in $L^2$ by the inequalities
\begin{align*}
\|(\A-\A')(\nabla-i(\A+\A')/2)v'\|_2 &\lesssim \|\A-\A'\|_4 \|(\nabla-i(\A+\A')/2)v';H^{3/4} \| \\
&\lesssim \|\A-\A'\|_4 \langle \|\A+\A';H^1\|\rangle\|v';H^{7/4}\|, \\
\|\omega^{-2}(\textrm{Re}(u-u')(\overline{u+u'}))v'\|_2&\lesssim \|u-u'\|_2\|u+u';H^1\|\|v';H^1\|,
\end{align*}
where we have used \eqref{nablaes1} and Lemma \ref{phies}.
Therefore, by the H\"older inequality for the time variable we have 
\begin{equation}\label{L2es}
\|v-v';L^\infty L^2\|\lesssim \|u_0-u_0'\|_2 
+T^{3/4}\langle \|v';L^\infty H^{7/4}\|\vee l_S\vee l_M % \|u;L^\infty H^1\|\vee \|\A;L^\infty H^{3/4}\| \vee \|\A';L^\infty H^{3/4}\|
\rangle^2 d(u,\A, u', \A'). %(\|\A-\A';L^4L^4\|\vee \|u-u';L^\infty L^2\|).
\end{equation}
On the other hand, by Lemma \ref{wavest} we have 
\begin{align*}
\|\B-\B';L^\infty H^{1/2}\cap L^4L^4\|&\lesssim \|(\A_0-\A_0',\A_1-\A_1');H^{1/2}\oplus H^{-1/2}\|\\
&\quad  +\|\A-\A';L^1H^{-1/2}\| +\|P(\J-\J');L^{4/3}L^{4/3}\|.
\end{align*}
For the last term we have the identity
\begin{equation}
P(\J-\J')=2P\Im\bigl\{(\overline{u-u'})(\nabla -i\A)u-i\bar{u}'(\A-\A')u
-(u-u')\overline{(\nabla -i\A') u'}\bigr\}
\label{difJ}
\end{equation}
since $P\{\bar{u'}(\nabla-i\A)(u-u')\}=-P\{(u-u')\overline{(\nabla-i\A)u'}\}$.
We apply the H\"older and the Sobolev inequalities to \eqref{difJ} together with \eqref{nablaes1}
and obtain
\[ \|P(\J-\J');L^{4/3}\| \lesssim \langle l_S \vee l_M \rangle^2 (\| u-u'\|_2 \vee \|\A-\A'\|_4). \]
Therefore the H\"older inequality for the time variable yields
%the  following estimates:
%\begin{align*}
%\|(\overline{u-u'})(\nabla-i\A) u;L^{4/3}\| &\lesssim \|u-u'\|_2\|(\nabla-i\A) u ;H^{3/4}\| \\
%&\lesssim \|u-u'\|_2 \langle \|\A:H^1\| \rangle \|u;H^{7/4}\|, \\
%\|\bar{u}'(\A-\A')u; L^{4/3}\| &\lesssim \|u';H^{7/4}\|\|\A-\A'\|_4\|u;H^{7/4}\|.
%\end{align*}
%
%We write 
%\begin{equation}
%\begin{split}
%P(\J-\J')&=2P\Im\{(\overline{u-u'})\nabla u+\bar{u}'\nabla(u-u')-i(\overline{u-u'})\A u \\
%&\qquad -i\bar{u}'(\A-\A')u-i\bar{u}'\A'(u-u')\}
%\end{split}
%\label{difJ}
%\end{equation}
%and estimate the right-hand side by 
%\begin{align*}
%\|(\overline{u-u'})\nabla u;L^{4/3}\| &\lesssim \|u-u'\|_2\|\nabla u \|_4, \\
%\|(\overline{u-u'})\A u;L^{4/3}\| &\lesssim \|u-u'\|_2\|\A\|_4\|u\|_\infty, \\
%\|\bar{u}'(\A-\A')u; L^{4/3}\| &\lesssim \|u'\|_2\|\A-\A'\|_4\|u\|_\infty.
%\end{align*}
%For $\bar{u}'\nabla(u-u')$ 
%we  use the formula $P(\bar{u}'\nabla(u-u'))=-P(\nabla \bar{u}'(u-u'))$
%and estimate it by the first inequality above.
%By the embeddings $H^{7/4}\embedding H^{1,4}\cap L^\infty$, 
%$H^{3/4}\embedding L^4$, we have
\[ \|\A-\A';L^1H^{-1/2}\| +\|P(\J-\J');L^{4/3}L^{4/3}\|
\lesssim T^{1/2}\langle l_S\vee l_M \rangle^2 d(u,\A, u', \A'). \]
Thus we obtain the required estimate.
\ebox
\thmskip
\Proofof{Theorem {\rm \ref{thmsc}}\textup{.} Part {\rm 1}}
We prove the unique existence of the solution for $s\ge 7/4$.
The case $s<7/4$ and the continuous dependence on the data
will be proved in later sections.
We consider the complete metric space $(Z_{s,\sigma},d)$
and the map $\Phi$
(see Proposition \ref{prop:itself} and \eqref{eq:d}
for the definition).
By Proposition \ref{prop:itself}, 
we can choose $l_S,l_M,\dots, L^s_S,L^\sigma_M,\tilde{L}^\sigma_M,T$ so that 
$\Phi$ is a map from $(Z_{s,\sigma},d)$ to itself.
On the other hand by Proposition \ref{metric}
\[ d(\Phi(u,\A),\Phi(u',\A'))\le C T^{1/2}\langle l_S\vee l_M\rangle^2 d(u,\A,u',\A'). \]
If we take suitable $T$, then $\Phi$ becomes a contraction mapping on $(Z_{s,\sigma},d)$.
This yields the unique existence of the fixed point $(u,\A)$.
This is the unique solution stated in the theorem.
Precisely the following is yet to be checked.
First, we check $(u,\A,\partial_t \A) \in C(I;X^{s,\sigma})$. 
Indeed we have $(\A,\partial_t \A) \in C(I;H^\sigma \oplus H^{\sigma-1})$
by Lemma \ref{W}.
On the other hand by virtue of Lemma \ref{lem:LSHm} and Corollary \ref{cor:Dtv}, we have 
$u \in C(I;H^s)\cap C^1 (I;H^{s-2})$ since $v=u$ is a solution to \eqref{v}.
Next, we have to check the uniqueness; 
we have used slightly different spaces from $C(I;X^{s,\sigma})$.
To this end we recall Lemma \ref{W}. From this lemma, 
we automatically have $\A \in L^4L^4$ if $(u,\A)$ is a solution in the required class.
Therefore the contraction argument above implies the uniqueness.
\ebox

%%%%%%%%%%% The contraction argument for $5/3\le s<7/4}
\section{The contraction argument for $5/3\le s<7/4$}
In this section we consider the following more complicated metric to refine the result on unique existence of the solution.
We put
\begin{equation}
\begin{split}\label{eq:tild}
\tild(u,\A,u',\A')&\equiv \max_{j=0,1}\|\partial^j_t(u-u');L^\infty H^{s-1-2j}\| \\
&\quad \vee\|\A-\A';L^q H^{2-s,r}\cap L^2L^\infty \cap L^\infty H^1\|\vee \|\partial_t(\A-\A');L^\infty L^2\|.
\end{split}
\end{equation}
Here $(q,r)=(6/(2s-1),3/(2-s))$. The space $L^q H^{2-s,r}$ is removed if $s=2$.
\begin{proposition}
\label{metric2}
Let $5/3\le s\le 2${\rm ,} $I=[0,T]$ with $0<T\le 1${\rm .}
Let $(u,\A), (u',\A')\in \tilde{Z}_{s,4/3}$ and let $(v,\B)$ and $(v',\B')$ be the solutions to
\eqref{v}{\rm -}\eqref{B} and $\eqref{v}'${\rm -}$\eqref{B}'${\rm,} respectively{\rm .}
Then the estimate 
\begin{align*}
&\tild(v,\B,v',\B') \le C(L)\|(u_0-u_0',\A_0-\A_0',\A_1-\A_1');X^{s-1,1+1/r}\| \\
& \quad+C(L)(\|(u-u')(0);H^{s-3}\|\vee\|(\A-\A')(0)\|_2)\|v';L^\infty H^s\| \\
& \quad+C(L)T^{1/2} \{\langle \|v';L^\infty H^s\|\vee \|\partial_tv';L^\infty H^{s-2}\|\rangle 
   \tild (u,\A,u',\A')+\|v-v';L^\infty H^{s-1}\|\}
\end{align*}
holds{\rm ,} where $L=l_S \vee l_M${\rm .}
If $s=2${\rm ,} the space 
$X^{s-1,1+1/r}$ in the estimate above is replaced by $X^{1,1+\delta}$ 
for sufficiently small $\delta>0${\rm .}
\end{proposition}

To prove this proposition, we need the following lemma which allows us to exchange
$\| v-v';L^\infty H^{s-1}\|$ and $\|\partial_t(v-v');L^\infty H^{s-3}\|$.
We state the lemma in general form although we use only \eqref{eq:sim2} to prove Proposition \ref{metric2}.
This is because we need \eqref{eq:sim1} to prove the continuous dependence of the solution on the data in Section 7.

\begin{lemma}\label{lem:sim}
Let $s>1/2${\rm ,} $\sigma \ge \max\{1,(2s-1)/4,s-2\}$ with $(s,\sigma)\neq (7/2,3/2)${\rm .}
Let $v$ and $v'$ be the solutions to $\eqref{v}$ and $\eqref{v}'${\rm,} respectively{\rm .}
Then the following inequality holds{\rm :}
\begin{align}
&\|v-v';H^s\| \notag
\\ 
&\quad\lesssim \|\partial_t(v-v');H^{s-2}\| +C(l)\{\|v-v'\|_2 +\| v';H^s\|(\|u-u';H^{s-2}\|+\|\A-\A';H^\sigma \|)\} \notag \\ 
&\quad\lesssim C(l)\{ \|v-v';H^s\|+\| v';H^s\|(\|u-u';H^{s-2}\|+\|\A-\A';H^\sigma \|)\}. \label{eq:sim1}
\end{align}
Here $l=\| u;H^s \|\vee\| u';H^s \|\vee\| \A;H^\sigma \|\vee\| \A';H^\sigma \|${\rm .}
Moreover{\rm ,} the following inequality holds for $3/2< s\le 2${\rm :}
\begin{align}
&\|v-v';H^{s-1}\| \notag \\
&\quad\lesssim \|\partial_t(v-v');H^{s-3}\| +C(l)\{\|v-v'\|_2 +\| v';H^s\|(\|u-u';H^{s-3}\|+\|\A-\A' \|_2)\} \notag \\
&\quad\lesssim C(l)\{\|v-v';H^{s-1}\|+\| v';H^s\|(\|u-u';H^{s-3}\|+\|\A-\A' \|_2)\}. \label{eq:sim2}
\end{align}
\end{lemma}
\Proof
In the beginning, we prove the first inequality of \eqref{eq:sim1}.
Applying \eqref{equi2} to \eqref{1difv}, we obtain  
\begin{align*}
\|\partial_t(v-v');H^{s-2}\| &\lesssim \|v-v';H^s\|+C(l)\|v-v'\|_2 
+\|(\calh+\phi-\calh'-\phi')v';H^{s-2}\|.
\end{align*}
For the estimate of $(\calh+\phi-\calh'-\phi')v'$, we need the following inequalities:
\begin{align}
\|(\A-\A')(\nabla-i(\A+\A')/2)v';H^{s-2}\|&\lesssim \|\A-\A';H^\sigma\|\langle\|\A+\A';H^\sigma\|\rangle\|v';H^s\|,
\label{Adifnv}\\
\|(\phi-\phi')v';H^{s-2}\|&\lesssim \|u-u';H^{s-2}\|\|u+u';H^s\| \|v';H^s\|.
\label{pdifv}
\end{align}
We can show \eqref{Adifnv} by the same method in the proof of Lemma \ref{equinorm}.
To prove \eqref{pdifv} for $s\ge2$, we use Lemma \ref{phies}.
If $1/2<s<2$, we also use duality.
%need
%the duality estimate as follows:
%\begin{align}
%|\langle \omega^{-2}((u-u')(\overline{u+u'})) v',\psi \rangle | 
%&=|\langle (u-u'),\omega^{-2}(\psi\bar{v}')(u+u') \rangle | \notag \\
%&\lesssim\|u-u';H^{s-2}\|\|\psi;H^{2-s}\| \|v';H^s\|\|u+u';H^s\|. 
%\end{align}
Therefore the first inequality of \eqref{eq:sim1} has been established.
The second inequality of  \eqref{eq:sim1} can be proved analogously.
To prove \eqref{eq:sim2}, we use similar estimates with $s$ replaced by $s-1$,
but we also need to use the estimate
\[ \|(\A-\A')(\nabla-i(\A+\A')/2)v';H^{s-3}\|
\lesssim \|\A-\A' \|_2\langle\|\A+\A';H^1\|\rangle\|v';H^s\| \]
instead of \eqref{Adifnv}.
This inequality is proved by the embedding $L^{6/(9-2s)}\embedding H^{s-3}$ and \eqref{nablaes1}.
\ebox
\thmskip

\Proofof{Proposition \rm{\ref{metric2}}}
Taking the difference of \eqref{v2} and the corresponding equation for $v'$, we have
\begin{equation}\label{2diffv}
\begin{split}
i\partial_t^2(v-v')&= (\calh+\phi)\partial_t(v-v')+(\calh+\phi-\calh'-\phi')\partial_tv' \\
&\quad+(2i\partial_t\A(\nabla-i\A)+\partial_t\phi)(v-v') \\
&\quad+(2i\partial_t\A(\nabla-i\A)+\partial_t\phi-2i\partial_t\A'(\nabla-i\A')-\partial_t\phi')v' \\
&\equiv (\calh+\phi)\partial_t(v-v')+\sum_{j=1}^3 f_j. 
\end{split}
\end{equation}
We convert this equation to the integral form by Lemma \ref{lem:inhomogeneous}
regarding $f\equiv \sum_{j=1}^3 f_j$ as the inhomogeneous term,
and take the $H^{s-3}$-norm of the both-sides. 
%
%Then 
%\begin{equation}
%\partial_t(v-v')(t)=U(t,0)\partial_t(v-v')(0)-i\int_0^tU(t,\tau)f(\tau)d\tau,
%\label{1difint}
%\end{equation}
%where $f\equiv \sum_{j=1}^3 f_j$.
Then we obtain
\begin{equation}\label{eq:intdiffv}
\|\partial_t(v-v')(t);H^{s-3}\|\le K_{s-3}\{\|\partial_t(v-v')(0);H^{s-3}\|+
\int_0^t\|f(\tau);H^{s-3}\|d\tau\}.
\end{equation}
$K_{s-3}\le C(L)$ clearly follows from Lemma \ref{lem:LSHm}. We begin with the estimate of the inhomogeneous term.
We have 
\[
\|f_1;H^{s-3}\| \lesssim \langle L\rangle\|\partial_tv';H^{s-2}\|(\|\A-\A';H^{s-2,r}\cap L^\infty\|\vee\|u-u';H^{s-1}\|).
\]
If $s=2$, we do not need $H^{2-s,r}$.
We can show this inequality by the duality argument such as 
\begin{align*}
&|\langle(\A-\A')(\nabla-i(\A+\A')/2)\partial_tv',\psi\rangle| \\
&\quad\le \|\partial_tv';H^{s-2}\|\|(\A-\A')(\nabla-i(\A+\A')/2)\psi;H^{2-s}\| \\
&\quad\lesssim\|\partial_tv';H^{s-2}\|\|\A-\A';H^{2-s,r}\cap L^\infty\|\|(\nabla-i(\A+\A')/2)\psi;H^{2-s}\|\\
&\quad\lesssim \|\partial_tv';H^{s-2}\|  \|\A-\A';H^{2-s,r}\cap L^\infty\| \langle\|\A+\A';H^1\|\rangle \|\psi;H^{3-s}\|
\end{align*}
and
\begin{align*}
|\langle (\phi-\phi')\partial_tv',\psi\rangle| 
&\le \|\partial_tv';H^{s-2}\|\|({\phi-\phi'})\psi;H^{2-s}\| \\
&\lesssim \|\partial_tv';H^{s-2}\|\|u-u';H^{s-1}\| \|u+u'\|_2 \|\psi;H^{3-s}\|,
\end{align*}
where we have used \eqref{nablaes1} and Lemma \ref{phies}. 
We next estimate $f_2$. We have 
\[
\|f_2;H^{s-3}\|\lesssim \|v-v';H^{s-1}\|\langle \|\partial_t\A\|_3\rangle\langle L\rangle^2
\]
again by the duality argument as follows:
\begin{align*}
|\langle \partial_t\A(\nabla-i\A)(v-v'),\psi\rangle| &\le \|\partial_t\A(v-v');H^{s-2}\|\|(\nabla-i\A)\psi;H^{2-s}\| \\
&\lesssim \|\partial_t\A\|_3\|v-v';H^{s-1}\|\|\psi;H^{3-s}\|\langle \|\A;H^1\|\rangle, \\
|\langle \omega^{-2}(\partial_t u \Bar{u})(v-v'),\psi\rangle| 
&\le \|\partial_t u;H^{s-2}\|\|\omega^{-2}((\overline{v-v'})\psi)u;H^{2-s}\| \\
&\lesssim \|\partial_tu;H^{s-2}\|\|v-v';H^{s-1}\|\|\psi;H^{3-s}\|\|u;H^1\|.
\end{align*}
For $f_3$, we have 
\begin{align*}
\|f_3;H^{s-3}\| \lesssim \langle L\rangle \|v';H^s\| \max_{j=0,1} 
\{\|\partial_t^j(\A-\A');H^{1-j}\|\vee \|\partial_t^j(u-u');H^{s-1-2j}\| \}  
\end{align*}
by the following estimates:
\begin{align*}
\|\partial_t\A(\A-\A')v';H^{s-3}\| &\lesssim \|\partial_t\A\|_2\|\A-\A'\|_6\|v';H^{s-1}\|, \\
\|\partial_t(\A-\A')(\nabla-i\A')v';H^{s-3}\| &\lesssim \|\partial_t(\A-\A')\|_2\|(\nabla-i\A')v';H^{s-3/2}\|\\
&\lesssim\|\partial_t(\A-\A')\|_2\|v';H^{s-1/2}\|\langle\|\A';H^1\|\rangle, \\
|\langle\omega^{-2}(\partial_t(u-u')\Bar{u})v',\psi\rangle| 
&\le \|\partial_t(u-u');H^{s-3}\|\|\omega^{-2}(\bar{v}'\psi) u;H^{3-s}\| \\
&\lesssim \|\partial_t(u-u');H^{s-3}\|\|v';H^s\|\|\psi;H^{3-s}\|\|u;H^s\|, \\
|\langle\omega^{-2}(\partial_t u' (\overline{u-u'}))v',\psi\rangle| 
&\le \|\partial_t u';H^{s-2}\|\|\omega^{-2}(\bar{v'}\psi)(u-u');H^{2-s}\| \\
&\lesssim \|\partial_tu';H^{s-2}\|\|v';H^s\|\|\psi;H^{3-s}\|\|u-u';H^{s-1}\|.
\end{align*}
Therefore we obtain 
\begin{align*}
\|F;L^1 H^{s-3}\| &\le 
C(L)T^{1/2}(\|v';L^\infty H^s\|\vee\|\partial_tv';L^\infty H^{s-2}\|)\tild(u,\A,u',\A') \\
&\quad+C(L)T^{5/6}\|v-v';L^\infty H^{s-1}\|.
\end{align*}
The estimate for $\| \partial_t(u-u');L^\infty H^{s-3}\|$ is completed
if we apply \eqref{eq:sim2} to the right-hand side of \eqref{eq:intdiffv}.
We also need the estimate for $\| u-u';L^\infty H^{s-1}\|$.
By virtue of \eqref{eq:sim2}, we only have to estimate
$\|u-u';L^\infty H^{s-3}\|$, $\|\A-\A' ;L^\infty L^2\|$ and $\| v-v';L^\infty L^2\|$.
To estimate the first two norms, we use the trivial inequality
$\|u-u';H^{s-3}\|\le\|(u-u')(0);H^{s-3}\|+\|\partial_t(u-u');L^1 H^{s-3}\|$
and the corresponding one for $\A-\A'$.
The norm $\| v-v';L^\infty L^2\|$ is estimated by using \eqref{difv},
but the estimate of the second term of the right-hand side is slightly different.
In this case we use the inequality
\begin{align*}
\| (\calh+\phi-\calh'-\phi')v' \|_2 &\lesssim \|\A-\A'\|_\infty\|v';H^1\|\langle\|\A+\A';H^1\|\rangle \\
&\quad+\|u-u';H^{1/2}\|\|u+u';H^{1}\|\|v'\|_2
\end{align*}
obtained by \eqref{nablaes1} and Lemma \ref{phies}. By the H\"older inequality for the time variable we have 
\[
\|v-v';L^\infty L^2\| \lesssim \|u_0-u_0'\|_2+T^{1/2}\tild(u,\A,u',\A')\|v';H^1\|\langle L\rangle.\]
Collecting the estimates above, we obtain
\begin{align*}
&\|v-v';L^\infty H^{s-1}\| \\
&\quad\le C(L)\bigl[
\|u_0-u_0';H^{s-1}\|+(\|(u-u')(0);H^{s-3}\|\vee\|(\A-\A')(0)\|_2)\|v';L^\infty H^s\| \\
&\qquad +T^{1/2}\{\langle\|v';L^\infty H^s\|\vee\|\partial_tv';L^\infty H^{s-2}\|\rangle\tild(u,\A,u',\A')
+\|v-v';L^\infty H^{s-1}\|\}\bigr].
\end{align*}
The estimate for the Schr\"odinger part has been completed.
We proceed to the Maxwell part. 
First we consider the case $s<2$; later we mention how to modify the proof when $s=2$.
We begin with the estimate in $L^qH^{2-s,r}$ with $s<2$.
We put $(\tilq,\tilr)\equiv (6/(4s-5),3/(4-2s))$, which is an admissible pair.
By Lemma \ref{wavest}, we have 
\begin{align*}
\|\B-\B'; L^q H^{2-s,r}\| 
&\lesssim \|(\A_0-\A_0',\A_1-\A_1');H^{\sigma_1}\oplus H^{\sigma_1-1}\| \\
&\quad+T\|\A-\A';L^\infty H^1\|+\|P(\J-\J');L^{\tilq'} H^{\sigma_1-1+\beta(\tilr),\tilr'}\|, 
\end{align*}
where $\sigma_1=1+1/r$ and $(\tilq,\tilr)$ is an admissible pair.
We should estimate the last term in the right-hand side.
We decompose $P(\J-\J')$ as in \eqref{difJ}.
%We only have to estimate (essentially) $(\overline{u-u'})(\nabla-i\A)u$ and $\bar{u}'(\A-\A')u$
%by virtue of the property $P(\bar{v}(\nabla-i\A)u)=-P(u\overline{(\nabla-i\A)v})$.
Since $s-1<\sigma_1<s$, we can choose $\tilr$ so that $\sigma_1-1+\beta(\tilr)=s-1$. 
With this choice, we have
\begin{equation}
\begin{split}
\|(\overline{u-u'})(\nabla-i\A)u;H^{s-1,\tilr'}\| 
&\lesssim \|u-u';H^{s-1}\|\|(\nabla-i\A)u;H^{s-1}\| \\
&\lesssim \|u-u';H^{s-1}\|\|u;H^s\|\langle\|\A;H^\sigma\|\rangle
\end{split}\label{unablau1}
\end{equation}
if $1/\tilr \le (s-1)/3$. This inequality holds provided $s\ge 5/3$.
We also have
\begin{equation}
\begin{split}
\|\bar{u}'(\A-\A')u;H^{s-1,\tilr'}\| 
&\lesssim \|u';H^{s-1,3}\|\|\A-\A'\|_6 \|u;1/2-1/\tilr\|  \\
&\quad +\|u'\|_\infty\ \|\A-\A';H^{s-1}\|\|u;1/2-1/\tilr \| \\
&\lesssim \|u';H^s\|\|u;H^s\|\|\A-\A';H^1\|.
\end{split}\label{unablau2}
\end{equation}
%\begin{align*}
%\|\B-\B'; L^q H^{2-s,r}\| 
%&\lesssim \|(\A_0-\A_0',\A_1-\A_1');H^{1+1/r}\oplus H^{1/r}\| \\
%&\quad+\|P(\J-\J');L^{\tilq'} H^{s-1,\tilr'}\|+T\|\A-\A';L^\infty H^1\|, 
%\end{align*}
%where we have used $H^1\embedding H^{1-s+\beta(r)}$ for the last term.
%If $5/3\le s <2$, we have
%\begin{equation}
%\begin{split}
%\|(\overline{u-u'})(\nabla-i\A)u;H^{s-1,\tilr'}\| &\lesssim \|u-u';H^{s-1}\|\|(\nabla-i\A)u;H^{s-1}\| \\
%&\lesssim \|u-u';H^{s-1}\|\|u;H^s\|\langle\|\A;H^1\|\rangle
%\end{split}\label{unablau1}
%\end{equation}
%{and}
%\begin{equation}
%\begin{split}
%\|\bar{u}'(\A-\A')u;H^{s-1,\tilr'}\| &\lesssim \|u';H^{s-1,6}\|\|\A-\A'\|_3\|u;1/2-1/\tilr \| \\
%&\quad +\|u'\|_\infty\ \|\A-\A';H^{s-1}\|\|u;1/2-1/\tilr \| \\
%&\lesssim \|u';H^s\|\|u;H^s\|\|\A-\A';H^1\|.
%\end{split}\label{unablau2}
%\end{equation}
Using these estimates together with the H\"older inequality for the time variable, we have
%\[
%\|P(\J-\J');L^{\tilq'} H^{s-1,\tilr'}\|\le T^{1/2}C(L)\tild(u,\A,u',\A').
%\]
%Therefore 
\begin{equation}
\begin{split}\label{eq:Bdifest}
\|\B-\B'; L^q H^{2-s,r}\| 
&\lesssim \|(\A_0-\A_0',\A_1-\A_1');H^{1+1/r}\oplus H^{1/r}\| \\
&\quad+T^{1/2}C(L)\tild(u,\A,u',\A').
\end{split}
\end{equation}
We can estimate $\| \B-\B';L^\infty H^1 \| \vee \| \partial_t(\B-\B');L^\infty L^2 \|$
in the same way. Indeed
these norms are estimated by the right-hand side of \eqref{eq:Bdifest}
with $H^{1+1/r}\oplus H^{1/r}$ replaced by $H^1\oplus L^2$.
Next we mention the estimate in $L^2L^\infty$;
the exponent $(2,\infty)$ is the prohibited endpoint.
However, for any admissible pair $(q_0, r_0)$ and any $\varepsilon>0$, we have 
\[ \|\B-\B';L^2L^\infty\|\lesssim \|\B-\B';L^{q_0}H^{3/r_0+\varepsilon, r_0}\|. \]
Here we have used the Sobolev inequality for the spatial variable and
the H\"older inequality for the time variable together with $0<T\le 1$.
Therefore we have by Lemma \ref{wavest}
\begin{equation}\label{eq:L2Linf}
\begin{split}
\|\B-\B';L^2L^\infty\|
&\lesssim \|(\A_0-\A_0',\A_1-\A_1');H^{\tilde{\sigma}_1}\oplus H^{\tilde{\sigma}_1-1}\|\\
&\quad +\|P(\J-\J');L^{\tilq_0'} H^{\tilde{\sigma}_1-1+\beta(\tilr_0),\tilr_0'}\|+T\|\A-\A';L^\infty H^1\|,
\end{split}
\end{equation}
where $\tilde{\sigma}_1=1+\varepsilon+1/r_0$ and $(\tilq_0,\tilr_0)$ is an admissible pair. 
If $8/5<s\le 2$, we can take $\varepsilon>0$ and $r_0,\tilr_0 \in [2,\infty)$ satisfying
\[ 0<1/r_0\le \min\{1/r-\varepsilon,(5s-8)/3-\varepsilon\},\quad \tilde{\sigma}_1-1+\beta(\tilr_0)=s-1. \]
With such a choice, we have 
$H^{1+1/r}\oplus H^{1/r} \embedding H^{1+\varepsilon+1/r_0}\oplus H^{\varepsilon+1/r_0}$
and 
\[
\|P(\J-\J');L^{\tilq'_0}H^{s-1,\tilr'_0}\| \lesssim T^{1/2}C(L)\tild(u,\A,u',\A'),
\]
since the estimates \eqref{unablau1} and \eqref{unablau2} hold with $\tilr$ replaced by $\tilr_0$ in the same manner.
Thus we have obtained the required result for the Maxwell part when $s<2$.
For $s=2$, we do not need $L^qH^{2-s,r}$. 
The estimate in $L^2L^\infty$ is almost same.
We only have to replace the condition $1/r_0\le 1/r-\varepsilon$ by $1/r_0\le \delta-\varepsilon$.
The norms $\|\B-\B';L^\infty H^1 \|$ 
and $\| \partial_t(\B-\B');L^\infty L^2 \|$ are also bounded by the right-hand side of \eqref{eq:L2Linf}
for the sake of Lemma \ref{wavest}.
Collecting the above estimates, we obtain the 
required result.
\ebox
\thmskip

\Proofof{Theorem {\rm \ref{thmsc}}\textup{.} Part {\rm 2}}
Here we treat the case $5/3\le s<7/4$.
We put
\[ \mathcal{B}\equiv \{ (u,\A) \in \tilde{Z}_{s,\sigma}; (u(0),\A(0),\partial_t\A(0))=(u_0,\A_0,\A_1) \} \]
and consider the complete metric space $(\mathcal{B},\tilde{d})$.
By Propositions \ref{prop:itself5/3} and \ref{metric2}, $\Phi$ is a map from $\mathcal{B}$ to itself
and
\[
\tild(v,\B,v',\B') \le C(l_M,l_S) T^{1/2}(\tild (u,\A,u',\A')+\|v-v';L^\infty H^{s-1}\|)
\]
with $(v,\B)=\Phi(u,\A)$ and $(v',\B')=\Phi(u',\A')$.
Therefore $\Phi$ is a contraction mapping if we take sufficiently small $T$.
\ebox

\begin{remark}
This proof is still valid if $7/4\le s\le 2$. Therefore the solution obtained in Section 5 
actually belongs to $\tilde{Z}_{s,\sigma}$ for some $l_S, l_M$.
\end{remark}

%%%%%%%%%%%%%%%%  Continuous dependence on the data
\section{Continuous dependence on the data}
In this section we prove the continuous dependence of the solution on the data,
which is usually the most delicate part of the theory of well-posedness.
Our method is essentially based on \cite{KP87}.
For a while, we assume the following:
\begin{assumption}\label{A2}\rm
(1) $(s,\sigma)$ satisfies 
$5/3\le s<\infty$, $(s,\sigma)\neq (5/2,7/2), (7/2,3/2)$ and 
\[ \max\{4/3,s-2,(2s-1)/4\} \le \sigma\le \min\{s+1, (5s-2)/3\}; \]

(2)
$I=[0,T]$ with $0<T\le 1$;

(3)
$(u,\A),(u',\A')$ are the solutions of \textbf{MS-C} with data 
$(u_0,\A_0,\A_1)$, $(u_0',\A_0',\A_1')\in X^{s,\sigma}$, respectively;

(4)
There exists a positive constant $L$ such that  $(u,\A)$ in (3) satisfies
\[\max_{j=0,1}(\| \partial_t^j u ;L^\infty H^{s-2j}\|\vee\|\partial_t^j\A;L^\infty H^{\sigma-j}\|)
\vee\| \partial_t\A;L^6L^3\|\le L. \]
$(u,\A)$ also satisfies the estimate 
$\| \partial_t^2 u ;L^\infty H^{s-4}\|\vee\|\partial_t^2\A;L^\infty H^{\sigma-2}\|\le L$ if $s>2$.
Moreover, $(u',\A')$ in (3) 
satisfies the estimates above with $(u,\A)$ replaced by $(u',\A')$.
\end{assumption}

To prove the continuous dependence, it is sufficient to show that
\begin{align*}
&D^{s,\sigma}(u,\A,u',\A') \\
&\quad\equiv \max_{j=0,1}\{\|\partial_t^j(\A-\A');L^\infty H^{\sigma-j}\cap L^6H^{\sigma-1/2-j,3}\| 
\vee \|\partial_t^j(u-u');L^\infty H^{s-2j}\|\}
\end{align*}
converges to 0 when $(u_0',\A_0',\A_1')$ tends to $(u_0,\A_0,\A_1)$ in $X^{s,\sigma}$.
To this end we also need for $s\le 2$
\begin{align*}
E^s(u,\A,u',\A')\equiv 
\|\partial_tu';L^\infty H^{s-1}\|(\|\A-\A';L^qH^{2-s,r}\cap L^2L^\infty\|\vee \|u-u';L^\infty H^{s-1}\|),
\end{align*}
where $(q,r)=(6/(2s-1),3/(2-s))$. The pair $(q,r)$ is admissible.
The space $L^qH^{2-s,r}$ in $E^s$ is removed if $s=2$.

%%%%%% Lemma {cddatasle2}
\begin{lemma}
\label{cddatasle2}
We assume {\rm \ref{A2}} with $s\le 2${\rm .}
%Moreover let $(u,\A),(u',\A')\in \tilde{Z}_{s,\sigma}${\rm .}
Then the estimate 
\begin{align*}
D&\le C(L)\{\|(u_0-u_0',\A_0-\A_0',\A_1-\A_1');X^{s,\sigma}\| \\
&\quad+T^{1/2}(E+D)+\|u-u';L^\infty L^2\|\}
\end{align*}
holds{\rm ,} where $D=D^{s,\sigma}(u,\A,u',\A'), E=E^s(u,\A,u',\A')${\rm .}
\end{lemma}
\Proof
The required result is obtained analogously to the proof of Proposition \ref{metric2}.
Therefore we again begin with the expression \eqref{2diffv} with $u=v, u'=v'$.
In the following, we can remove $H^{2-s,r}$ and $L^qH^{2-s,r}$ if $s=2$.
First, we estimate $f=\sum_{j=1}^3 f_j$ in the inhomogeneous term. We have
\[ \| f_1;H^{s-2}\| \le C(L) \| \partial_t u';H^{s-1}\| \|\A-\A';H^{2-s,r}\cap L^\infty\| 
+C(L)\|u-u';H^s\| \]
by virtue of
\begin{align*}
&|\langle(\A-\A')(\nabla-i(\A+\A')/2)\partial_tu',\psi\rangle| \\
&\quad\le \|(\nabla-i(\A+\A')/2)\partial_tu';H^{s-2}\|\|(\A-\A')\psi;H^{2-s}\|\\
&\quad\lesssim \|\partial_tu';H^{s-1}\|\langle\|\A+\A';H^1\|\rangle\|\A-\A';H^{2-s,r}\cap L^\infty\|\|\psi;H^{2-s}\|, 
\\
&|\langle (\phi-\phi')\partial_tu',\psi\rangle| \lesssim \|u-u';H^s\|(\|u;H^s\|+\|u';H^s\|)\|\partial_tu';H^{s-2}\|\|\psi;H^{2-s}\|.
%&\|(\phi-\phi')\partial_tu';H^{s-2}\| \lesssim \|u-u';H^s\|\|u+u';H^s\|\|\partial_tu';H^{s-2}\|.
\end{align*}
We have
\[ 
\| f_2;H^{s-2}\| \le C(L) \|\partial_t \A \|_3 \|u-u';H^s\|
\]
by
\begin{align*}
\|\partial_t \A(\nabla-i\A)(u-u');H^{s-2}\|
&\lesssim \|\partial_t\A\|_3\|(\nabla-i\A)(u-u');H^{s-1}\| \\
&\lesssim \|\partial_t\A\|_3\|u-u';H^s\|\langle\|\A;H^1\|\rangle ,\\
|\langle\omega^{-2}(\partial_tu\bar{u})(u-u'),\psi\rangle|
&\lesssim \|\partial_tu;H^{s-2}\|\|u-u';H^s\|\|\psi;H^{2-s}\|\|u;H^s\|.
\end{align*}
We have
\[ \| f_3;H^{s-2}\| \le C(L) \{\|\A-\A';H^1\| 
+\|\partial_t(\A-\A')\|_3+\max_{j=0,1}\|\partial_t^j(u-u');H^{s-2j}\|\} \]
by
\begin{align}
\|\partial_t\A'(\A-\A')u';H^{s-2}\|
&\lesssim \|\partial_t\A'\|_2\|\A-\A';H^1\| \|u';H^s\|, \notag\\
\|\partial_t(\A-\A')(\nabla-i\A)u';H^{s-2}\|
&\lesssim \|\partial_t(\A-\A')\|_3\|u';H^s\|\langle\|\A;H^1\|\rangle, \label{eq:uselater}\\
|\langle\omega^{-2}(\partial_t(u-u')\bar{u}))u',\psi\rangle|
&\lesssim \|\partial_t(u-u');H^{s-2}\|\|u;H^s\|\|u';H^s\|\|\psi;H^{2-s}\|, \notag\\
|\langle\omega^{-2}(\partial_t u'(\overline{u-u'}))u',\psi\rangle|
&\lesssim \|\partial_t u';H^{s-2}\|\|u-u';H^s\|\|u';H^s\|\|\psi;H^{2-s}\|. \notag
\end{align}
Therefore by these inequalities together with the H\"older inequality for the time variable,
%
%$\| f;L^1 H^{s-2}\|$ is bounded by
%\begin{align*}
%\|\partial_t(u-u');H^{s-2}\| \lesssim C(L)\|\partial_t(u-u')(0);H^{s-2}\|+T^{1/2}C(L)(E+D)
%&
%T^{1/2}C(L)\bigl\{
%\|\partial_t u';L^\infty H^{s-1}\|\|\A-\A';L^qH^{2-s,r}\cap L^2L^\infty\| \\
%& \qquad+\|\A-\A';L^\infty H^1\|\vee \|\partial_t(\A-\A');L^6L^3\| \vee 
%\max_{j=0,1}\|\partial_t^j(u-u');L^\infty H^{s-2j}\|\bigr\}.
%\end{align*}
%
\[
\|\partial_t(u-u');H^{s-2}\| \lesssim C(L)\|\partial_t(u-u')(0);H^{s-2}\|+T^{1/2}C(L)(E+D).
\]
We obtain the estimate for $\|\partial_t(u-u');L^\infty H^{s-2}\|$ by using \eqref{eq:sim1} to the first term of the right-hand side.
Next we estimate  $\|u-u';L^\infty H^s\|$. By  the interpolation inequality to \eqref{eq:sim1}, we have
\[ \| u-u';H^s\|\lesssim \|\partial_t(u-u');H^{s-2}\|+C(L)\|u-u'\|_2+C(L)\|\A-\A'; H^\sigma\|.\]
Therefore
\begin{align}
\max_{j=0,1}\|\partial_t^j (u-u');L^\infty H^{s-2j}\|
&\le C(L)\{ \| (u_0-u_0',\A_0-\A_0');H^s \oplus H^\sigma\| +T^{1/2}(E+D)\nonumber \\
&\quad+\|\A-\A';L^\infty H^\sigma\| +\|u-u';L^\infty L^2\| \}. \label{diffues}
\end{align}
On the other hand, by the analogous argument in the proof of Lemma \ref{W}, we have 
\begin{align}
&\max_{j=0,1}\|\partial_t^j(\A-\A');L^\infty H^{\sigma-j}\cap L^qH^{\sigma-j-\beta(r),r}\| \nonumber\\
&\quad\lesssim \|(\A_0-\A_0',\A_1-\A_1');H^\sigma\oplus H^{\sigma-1}\| \nonumber\\
&\qquad+C(L)T^{1/2}(\|\A-\A';L^\infty H^\sigma\|\vee \|u-u';L^\infty H^s\|)
%& \qquad \times \langle \|\A;L^\infty H^\sigma\|\vee \|\A';L^\infty H^\sigma\|\vee 
%\|u;L^\infty H^s\|\vee \|u';L^\infty H^s\|\rangle^2.
\label{diffAes}
\end{align}
%We replace $\|\A-\A';L^\infty H^\sigma\|$ in \eqref{diffues} by the right-hand of \eqref{diffAes}.
%After that we add \eqref{diffues} and \eqref{diffAes}.
Thus we have obtained the desired result.
\ebox
\thmskip

If $s>2$, instead of $D^{s,\sigma}$ and $E^s$ we need 
\begin{align*}
\tilD^{s,\sigma}(u,\A,u',\A')&\equiv
D(u,\A,u',\A')\vee \|\partial^2_t(\A-\A');L^\infty H^{\sigma-2}\| \vee \|\partial_t^2(u-u');L^\infty H^{s-4}\|, \\
\tilE^s(u,\A,u',\A')&\equiv 
\|\partial_t^2u';H^{s-3}\|(\|\A-\A';L^{q}H^{|s-3|,r}\cap L^2L^\infty\|\vee \|u-u';L^\infty H^{s-1}\|),
\end{align*}
where $1/r=1/2-(1/2-|s-3|/3)_+$. 
We choose $q$ so that $(q,r)$ is an admissible pair.
The space $L^qH^{|s-3|,r}$ in $\tilE^s$ is removed if $s=3$.
\begin{lemma}
\label{cddatasge2}
Let assume {\rm\ref{A2}} with $s>2${\rm .}
Then the  estimate 
\begin{align*}
 \tilD \lesssim C(L)\{\|(u_0-u_0',\A_0-\A_0',\A_1-\A_1');X^{s,\sigma}\| +T^{1/2}(\tilE+\tilD) %+\max_{j=0,1}\|\partial_t^j(u-u');L^\infty H^{s-1-2j}\|
+\|u-u';L^\infty L^2\|\}
\end{align*}
holds{\rm ,} where $\tilD=\tilD^{s,\sigma}(u,\A,u',\A')$ and $\tilE=\tilE^s(u,\A,u',\A')${\rm .}
%$\varepsilon>0$ is sufficiently small number dependent only on $s${\rm .}
\end{lemma}

\Proof
By taking the difference of \eqref{v3} with $v=u$ and the corresponding equation for $u'$, we have 
\begin{align*}
&i\partial_t^3(u-u') = (\calh+\phi)\partial_t^2(u-u')+(\calh+\phi-\calh'-\phi')\partial_t^2u'\\
&\qquad+2(2i\partial_t\A(\nabla-i\A)+\partial_t\phi)\partial_t(u-u')\\
&\qquad +2(2i\partial_t\A(\nabla-i\A)+\partial_t\phi-2i\partial_t\A'(\nabla-i\A')-\partial_t\phi')\partial_tu' \\
&\qquad+(2i\partial_t^2\A(\nabla-i\A)+2(\partial_t\A)^2+\partial_t^2\phi)(u-u') \\
&\qquad+(2i\partial_t^2\A(\nabla-i\A)+2(\partial_t\A)^2+\partial_t^2\phi-2i\partial_t^2\A'(\nabla-i\A')-2(\partial_t\A')^2-\partial_t^2\phi')u'\\
&\quad\equiv (\calh+\phi)\partial_t^2(u-u')+\sum_{j=1}^5 G_j.
\end{align*}
%We put $G\equiv \sum_{j=1}^5G_j$. 
%Then we have 
%\[
%\partial_t^2(u-u')(t)=U(t,0)\partial_t^2(u-u')(0)-i\int_0^tU(t,\tau)G(\tau)d\tau
%\]
%and therefore
Regarding $G\equiv \sum_{j=1}^5G_j$ as the inhomogeneous term,
we convert this equation to the integral form by virtue of Lemma \ref{lem:inhomogeneous},
and take the $H^{s-4}$-norm of the both-sides. Then we obtain
\begin{equation}
\|\partial_t^2(u-u');H^{s-4}\|\lesssim K_{s-4}\{\|\partial_t^2(u-u')(0);H^{s-4}\|+\int_0^t\|G(\tau);H^{s-4}\|d\tau\}.
\label{cdu2integral}
\end{equation}
We begin with the estimate of the inhomogeneous term. For $G_1$, we have 
\[
\|G_1;H^{s-4}\| \le C(L)\|\partial_t^2u';H^{s-3}\|(\|\A-\A';H^{|s-3|,r}\cap L^\infty\|\vee \|u-u';H^{s-1}\|),
\]
where $H^{|s-3|,r}$ is removed if $s=3$. 
To prove this, we should estimate 
\[ G_{1,1}=(\A-\A')(\nabla-i(\A+\A')/2)\partial_t^2u' \text{\quad and\quad} G_{1,2}=(\phi-\phi')\partial_t^2u'. \]
If $s>3$, we can rewrite $G_{1,1}=(\nabla-i(\A+\A')/2)(\A-\A')\partial_t^2u'$ 
since $\Div (\A-\A')=0$. After that, we use \eqref{nablaes1} and the Leibniz formula to estimate $G_{1,1}$.
If $2<s\le 3$,  we also use duality.
In both cases we have
\[ G_{1,1} \lesssim C(L)\|\partial_t^2u';H^{s-3}\| \|\A-\A';H^{|s-3|,r}\cap L^\infty\|. \]
We have $G_{1,2}\lesssim C(L)\| u-u'; H^{s-1}\|$.
This is obtained by Lemma \ref{phies}, together with duality if $s<3$.
%\begin{align*}
%|\langle (\A-\A')(\nabla-i(\A+\A')/2)\partial_t^2u',\psi\rangle| 
%&\lesssim \|\partial_t^2u';H^{s-3}\|\|(\A-\A')(\nabla-i(\A+\A')/2)\psi;H^{3-s}\|\\
%\|(\A-\A')(\nabla-i(\A+\A')/2)\psi;H^{3-s}\| 
%&\lesssim \|\A-\A';H^{3-s,r}\cap L^\infty\|\|\psi;H^{4-s}\|\langle\|\A+\A';H^1\|\rangle\\
%|\langle\omega^{-2}((u-u')\bar{u})\partial_t^2u',\psi\rangle|
%&\lesssim \|\partial_t^2u';H^{s-3}\|\|u-u';H^{s-1}\|\|u;H^{s-1}\|\|\psi;H^{4-s}\|
%\end{align*}
%for $2<s\le 3$, here $H^{3-s,r}$ is removed if $s=3$, and the estimates 
%\begin{align*}
%\|(\A-\A')(\nabla-i(\A+\A')/2)\partial_t^2u';H^{s-4}\| 
%&\lesssim \|(\A-\A')\partial_t^2u';H^{s-3}\|\langle\|\A+\A';H^\sigma\|\rangle \\
%\|(\A-\A')\partial_t^2u';H^{s-3}\|
%&\lesssim \|\A-\A';H^{s-3,r}\cap L^\infty\|\|\partial_t^2u';H^{s-3}\|
%\end{align*}
%for $s>3$.
For the estimates of $G_j$, $2\le j\le 5$, 
we can use the estimates for $F_j$, $1\le j\le 5$, in the proof of Lemma \ref{lem:HsestLS}.
Indeed we have
\begin{align*}
\sum_{j=2}^5\| G_j;H^{s-4}\| 
&\le C(L) \langle \| \partial_t \A\|_3 \rangle 
\max_{j=0,1}\{\| \partial_t^j (\A-\A');H^{\sigma-j} \| \vee \|\partial_t^j (u-u');H^{s-2j}\|\} \\
&\quad+C(L)\| \partial_t(\A-\A')\|_3+C(L)\| \partial_t^2(\A-\A');H^{\sigma-2} \|.
\end{align*}
Therefore we obtain
\[
\|G;L^1 H^{s-4}\| \le T^{1/2}C(L)\bigl\{ \tilD^{s,\sigma}+\tilE^s \bigr\}.
\]
To complete the estimate for the Schr\"odinger equation,
we need the following estimate:
\begin{align}
&\max_{j=0,1}\|\partial_t^j(u-u');H^{s-2j} \| \notag \\
&\quad\lesssim \|\partial_t^2 (u-u');H^{s-4} \|+C(L)\| u-u'\|_2
+C(L)\max_{j=0,1}\| \partial_t^j (\A-\A');H^{\sigma-j}\| \notag \\
&\quad\lesssim C(L)\{ \| u-u';H^s\| + \max_{j=0,1}\| \partial_t^j (\A-\A');H^{\sigma-j}\|\}. \label{sim3}
\end{align}
To this end we recall \eqref{2diffv}. 
Applying \eqref{equi2} and Lemma \ref{lem:sim} to this equation, we obtain
\begin{align*}
\max_{j=0,1}\|\partial_t^j(u-u');H^{s-2j} \|&\lesssim \|\partial_t^2 (u-u');H^{s-4} \|+C(L)\|u-u'\|_2 \\
&\quad+C(L)\|\A-\A';H^\sigma\|+\sum_{j=1}^3\|f_j;H^{s-4}\|, 
\end{align*}
where $f_j$, $j=1,2,3$, are defined in \eqref{2diffv} with $v,v'$ replaced by $u,u'$ respectively.
The estimate for $f_1$ is obtained by \eqref{equi1} and Lemma \ref{phies}.
%\[\| f_1;H^{s-4} \| \le C(L) (\| \A-\A';H^\sigma \| + \| u-u';H^{s-1}\|). \]
%We have $\| f_2;H^{s-4} \| \le C(L) \| u-u';H^{s-1}\|$
%\[ \| f_3;H^{s-4} \| \le C(L) 
%\max_{j=0,1} (\| \partial_t^j (\A-\A');H^{\sigma-j} \|\vee \| \partial_t^j (u-u');H^{s-2j-1} \|)
%\]
The estimate for $f_2$ and $f_3$ are obtained
similarly to the estimate for $F_6,F_7$ in the proof of Lemma \ref{lem:HsestLS}.
Indeed we have
\begin{equation} \label{eq:f1tof3}
\sum_{j=1}^3\|f_j;H^{s-4}\|\le 
C(L)\max_{j=0,1}(\|\partial_t^j(\A-\A');H^{\sigma-j}\|
\vee \|\partial_t^j(u-u');H^{s-1-2j}\|).
\end{equation}
The right-hand side of \eqref{eq:f1tof3} 
does not exceed 
\[ C(L)\max_{j=0,1}\|\partial_t^j(\A-\A');H^{\sigma-j}\| +C(L)\|u-u';H^{s-1}\| \]
again by Lemma \ref{lem:sim}.
Therefore using the interpolation inequality to $\|u-u';H^{s-1}\|$, we obtain the first inequality of \eqref{sim3}.
Similarly we can obtain the second.
Thus we have the following estimate for the Schr\"odinger part:
\begin{align*}
\max_{j=0,1,2} \|\partial_t^j(u-u');H^{s-2j}\| &\le C(L)\bigl[ \| (u_0-u_0',\A_0-\A_0',\A_1-\A_1');X^{s,\sigma}\| 
+T^{1/2}(\tilE+\tilD) \\
&\qquad+\max_{j=0,1} \| \partial_t^j (\A-\A');L^\infty H^{\sigma-j} \|
+ \| u-u';L^\infty L^2\|\}\bigr].
\end{align*}
The Maxwell part is easy to treat. Indeed, the estimate \eqref{diffAes} is still valid for $s>2$;
we also have $\| \partial_t^2 (\A-\A');L^\infty H^{\sigma-2}\| \le C(L) D^{s,\sigma}$ similarly as in Lemma \ref{W}.
Collecting the estimates both for the Schr\"odinger equation and the Maxwell, we obtain the desired result. 
\ebox

\begin{lemma}
\label{cddiffA}
Let assume {\rm \ref{A2}.}
Let $N(\A,\A')$ be defined by the following condition{\rm :}

{\rm (1)}
 If $5/3\le s<2${\rm ,} let $1/r\equiv(2-s)/3, \sigma_1\equiv 1+1/r${\rm ,} 
\[
N(\A,\A')\equiv\|\A-\A';L^\infty H^{\sigma_1}\cap L^qH^{2-s,r}\cap L^2L^\infty\|;
\]

{\rm (2)}
 If $s>2$ and $s\neq3${\rm ,} let $1/r\equiv1/2-(1/2-|s-3|/3)_+, \sigma_1\equiv|s-3|+\beta(r)${\rm ,}
\[
N(\A,\A')\equiv\|\A-\A';L^\infty H^{\sigma_1}\cap L^qH^{|s-3|,r}\cap L^2L^\infty\|;
\]

{\rm (3)} If $s=2$ or $s=3${\rm ,} let $\sigma_1\equiv 1+\varepsilon$ for sufficiently small $\varepsilon>0${\rm ,}
\[
N(\A,\A')\equiv\|\A-\A';L^\infty H^{\sigma_1}\cap L^2L^\infty\|.
\]
In any case we choose $q$ so that $(q,r)$ is an admissible pair{\rm .}
Then the following estimate holds{\rm :}
\begin{align*}
N(\A,\A')&\lesssim \|(\A_0-\A_0',\A_1-\A_1');H^{\sigma_1}\oplus H^{\sigma_1-1}\| \\
&\quad+C(L)T^{1/2}(\|u-u';L^\infty H^{s-1}\|\vee\|\A-\A';L^\infty H^{\sigma_1}\|).
\end{align*}
\end{lemma}

\Proof
Since $0<T\le 1$, there exists an admissible pair $(q_0,r_0)$ such that 
$L^{q_0}H^{3/r_0+0, r_0}\embedding L^2L^\infty$ and that $1+1/r_0<\sigma_1$.
Therefore by Lemma \ref{wavest}, we have 
\begin{align*}
N(\A,\A') &\lesssim \|(\A_0-\A_0',\A_1-\A_1');H^{\sigma_1}\oplus H^{\sigma_1-1}\|\\
&\quad +\|P(\J-\J');L^{\tilq}H^{\sigma_1-1+\beta(\tilr),\tilr'}\|+T\|\A-\A';L^\infty H^{\sigma_1-1}\| ,
\end{align*}
where $(\tilq,\tilr)$ is an admissible pair.
We estimate the middle term of the right-hand side as in the proof of Proposition \ref{metric2}.
We choose $\tilr$ so that $1/2-1/\tilr\ge (1/2-(s-1)/3)_+$ and that $\sigma_1-1+\beta(\tilr)\le s-1$.
Then  \eqref{unablau1} and \eqref{unablau2} still hold valid 
with $\|\A-\A';H^1\|$ in \eqref{unablau2} replaced by $\|\A-\A';H^{\sigma_1}\|$.
Therefore  we can obtain the desired result.
\ebox

\begin{lemma}
\label{cddata2s3}
Let assume {\rm\ref{A2}} with $2<s\le3${\rm .} 
Let $s_1=s-1${\rm ,} and $\sigma_1$ be defined in Lemma {\rm \ref{cddiffA}.} 
Then  the estimate 
\begin{equation}
D\le C(L)(\|(u_0-u_0',\A_0-\A_0',\A_1-\A_1');X^{s_1,\sigma_1}\|+T^{1/2}D)
\label{D2s3}
\end{equation}
holds{\rm ,} where $D\equiv D^{s_1,\sigma_1}(u,\A,u',\A')${\rm .}
%Especially for sufficiently small $T>0${\rm ,} the following estimate holds{\rm :}
%\begin{equation}
%\|u-u';L^\infty H^{s_1}\| \le C(L)\|(u_0-u_0',\A_0-\A_0',\A_1-\A_1');X^{s_1,\sigma_1} \|.
%\label{diffusm1}
%\end{equation}
\end{lemma}

\Proof
%The proof is the combination of Proposition \ref{metric}, Lemma \ref{cddiffA} 
%and a slight modification of the proof for Lemma \ref{cddatasle2}.
First, we prove the following inequality:
\begin{align}
D^{s_1,\sigma_1}&\lesssim \|(u_0-u_0',\A_0-\A_0',\A_1-\A_1');X^{s_1,\sigma_1}\| \notag\\
&\quad+C(L)T^{1/2}(E^{s_1}+D^{s_1,\sigma_1})+\| u-u';L^\infty L^2 \|. \label{eq:step1}
\end{align}
Here 
$E^{s_1}=
\|\partial_tu';L^\infty H^{s_1-1}\|(\|\A-\A';L^{q_1}H^{2-s_1,r_1}\cap L^2L^\infty\|\vee \|u-u';L^\infty H^{s_1-1}\|)
$  
and $(q_1,r_1)=(6/(2s-3),3/(3-s))$.
This is the same inequality as the assertion of Lemma \ref{cddatasle2}.
However we have assumed $s\ge 5/3$, $\sigma\ge 4/3$ in Lemma \ref{cddatasle2}; $(s_1,\sigma_1)$ does not satisfy these conditions. 
Therefore we have to modify the proof.
In the proof of \eqref{eq:step1},
we should distinguish the assumption for $(s,\sigma)$ and that for $(s_1,\sigma_1)$.
We do not need $s_1\ge 5/3$; 
we need $s\ge 5/3$ only to ensure the unique existence of the solution, 
and $s_1>1$ is sufficient to obtain \eqref{eq:step1}.
In the proof of Lemma \ref{cddatasle2}, the assumption $\sigma\ge 4/3$ is used to ensure the boundedness of 
$\| \partial_t\A;L^6L^3 \|$ and $\| \partial_t(\A-\A);L^6L^3 \|$
by virtue of the Strichartz estimate.
The former norm is still bounded because even in the present case $\sigma$ itself satisfies this condition.
The latter norm cannot be controlled if $\sigma_1<4/3$, but it appears only once in the estimate of $F_3$.
Therefore if we replace \eqref{eq:uselater} by the following inequality, we obtain the estimate for the Schr\"odinger part: 
\[ \|\partial_t(\A-\A')(\nabla-i\A)u';H^{s_1-2}\|
\lesssim \|\partial_t(\A-\A')\|_2\|u';H^s\|\langle\|\A;H^1\|\rangle. \]
The estimate for the Maxwell part is similar to that in the proof of Proposition \ref{metric2}.
Therefore \eqref{eq:step1} has been established.
On the other hand, we have
\begin{gather*}
\begin{split}
E^{s_1}&\le C(L) N(\A,\A')\vee\| u-u';L^\infty H^{s-1} \| \\ 
&\le C(L)\| \A_0-\A_0',\A_1-\A_1';H^{\sigma_1}\oplus H^{\sigma_1-1}\|
+C(L) D^{s_1,\sigma_1}, 
\end{split}\\
\| u-u';L^\infty L^2 \| \le \|u_0-u_0'\|_2+C(L)T^{1/2}d(u,\A,u',\A') 
\end{gather*}
by Lemma \ref{cddiffA} and by Proposition \ref{metric} respectively.
We obtain the desired result by collecting these estimates and using the inequality $\|\A-\A';L^4L^4\|\le N(\A,\A')$,
which is obtained by interpolation. 
\ebox

\begin{lemma}
\label{cddatas3}
Let assume {\rm \ref{A2}} with $s>3${\rm .}
Let $s_1=s-1$ and $(s_1,\sigma_\ast)$ satisfy $\sigma_\ast<\sigma$ and 
{\rm\ref{A2}-(1)} with $(s,\sigma)$ replaced by $(s_1,\sigma_\ast)${\rm .}
Then the estimate 
\[
\tilD\le C(L)(\|(u_0-u_0',\A_0-\A_0',\A_1-\A_1');X^{s_1,\sigma_\ast}\|+T^{1/2}\tilD)
\]
holds{\rm ,} 
where $\tilD\equiv \tilD^{s_1,\sigma_\ast}(u,\A,u',\A')${\rm .}
%Especially{\rm ,} for sufficiently small $T>0${\rm ,} the following estimate holds{\rm .}
%\[
%\|u-u';L^\infty H^{s_1}\| \le C(L)\|(u_0-u_0',\A_0-\A_0',\A_1-\A_1');X^{s_1,\sigma_\ast}\|
%\]
\end{lemma}

\Proof
By Lemma \ref{cddatasge2}, we have 
\begin{align*}
\tilD^{s_1,\sigma_\ast} &\le C(L)\{\|(u_0-u_0',\A_0-\A_0',\A_1-\A_1');X^{s_1,\sigma_\ast}\| \\
&\quad +T^{1/2}(\tilD^{s_1,\sigma_\ast}+\tilE^{s_1})
+\|u-u';L^\infty L^2\|\}.
\end{align*}
We have $\sigma_\ast >\sigma_1$ for $s>3$, where $\sigma_1$ is defined in Lemma \ref{cddiffA}.
Therefore we can prove the lemma similarly to Lemma \ref{cddata2s3}.
\ebox
\thmskip

\Proofof{Theorem {\rm\ref{thmsc}.} Part {\rm 3}}
%In Lemma \ref{lem:HsestLS}, if $s\ge5/2$ and $(u,\A)$ is the solution of \textbf{MS-C}, then we can have for small $\varepsilon>0$
%\begin{equation}
%\|\partial_t^2u;L^\infty H^{s-4}\|\lesssim K_{s-4}(\|u_0;H^s\|+T^{1/3} \max_{j=0,1,2}\|\partial_t^ju;L^\infty H^{s-2j}\|)\langle W_\sigma\vee R^0_{s-1}\vee R^1_{s-1}\rangle^\alpha.
%\label{cdd2u}
%\end{equation}
%Indeed, the above inequality follows from the modification in the proof of (2) for $F_5$ for $2<s<4$ such as 
%\begin{align*}
%|\langle \omega^{-2}(\partial_t^2u\bar{u})u, \psi\rangle|
%&\le |\langle \partial_t^2u, \omega^{-2}(\bar{u}\psi)u\rangle| \\
%&\lesssim \|\partial_t^2u;H^{s-4}\|\|u;H^{s-1}\|^2\|\psi;H^{4-s}\| \\
%\|\omega^{-2}(\partial_tu\partial_t\bar{u})u;H^{s-4}\|
%&\lesssim \|\partial_tu;H^{s-3}\|^2\|u;H^{s-1}\|\ \ \mbox{for }\ \ 3\le s<4 \\
%|\langle \omega^{-2}(\partial_tu\partial_t\bar{u})u,\psi\rangle|
%&\le |\langle \partial_tu, \omega^{-2}(\bar{u}\psi)\partial_tu\rangle| \\
%&\lesssim \|\partial_tu;H^{s-3}\|\|u;H^{s-1}\|\|\psi;H^{4-s}\|\|\partial_tu;H^{s-2}\| \ \ \mbox{for}\ \ 5/2\le s<3.
%\end{align*}
%By \eqref{remainder}, $\|u;L^\infty H^s\|$ is estimated by the right-hand side of \eqref{cdd2u}.
%By (3) of Lemma \ref{lem:HsestLS}, we obtain for sufficiently small $T>0$
%\[
%\max_{j=0,1,2}\|\partial_t^ju;L^\infty H^{s-2j}\|\le C(W_\sigma\vee \max_{j=0,1}R^j_{s-1})\|u_0;H^s\|.
%\]
Here we prove the continuous dependence of the solution on the data.
Let $(u,\A)$ and $(u',\A')$ be the solutions of \textbf{MS-C} with data
$(u_0,\A_0,\A_1)$ and $(u_0',\A_0',\A_1')$ respectively.
We consider the case that $(u_0',\A_0',\A_1')\to (u_0,\A_0,\A_1)$ in $X^{s,\sigma}$. Therefore we may assume 
that $(u_0',\A_0',\A_1')$ is bounded in $X^{s,\sigma}$,
accordingly Assumption \ref{A2} is satisfied for some $L>0$.
We put $\varepsilon = \sigma-\sigma_1$ if $s\le 3$, $\varepsilon = \sigma-\sigma_\ast$ if $s>3$.
$\sigma_1$ and $\sigma_\ast$ are given in Lemmas \ref{cddiffA} and \ref{cddatas3}.
In both cases $\varepsilon >0$.
In the following, we take $T>0$ sufficiently small.
We summarize the estimates used in the proof.
By Lemma \ref{lem:HsestLS} and Corollary \ref{cor:Dtv}, we have
\begin{equation}
\max_{j=0,1,2}\|\partial_t^ju;L^\infty H^{s+1-2j}\| \le C(L)\|u_0;H^{s+1}\|
\label{uHsplus1}
\end{equation}
and the corresponding estimate for $u'$.
Since $T>0$ is sufficiently small, we have by Proposition \ref{metric2} and Lemmas \ref{cddiffA}-\ref{cddatas3}
\begin{equation}
\|u-u';L^\infty H^{s-1}\|\vee N(\A,\A')\le C(L)\|(u_0-u_0',\A_0-\A_0',\A_1-\A_1');X^{s-1,\sigma-\varepsilon}\|.
\label{cduusminus1}
\end{equation}
By Lemmas \ref{cddatasle2} and \ref{cddatasge2} together with \eqref{uHsplus1} and \eqref{cduusminus1}, we also have 
\begin{align}
D^{s,\sigma}(u,\A,u',\A') &\le C(L)\{\|(u_0-u_0',\A_0-\A_0',\A_1-\A_1');X^{s,\sigma}\| \nonumber\\
&\quad+E^s(u,\A,u',\A')+\| u-u';L^\infty L^2\|\} \nonumber \\
&\le C(L)\{\|(u_0-u_0',\A_0-\A_0',\A_1-\A_1');X^{s,\sigma}\| \nonumber\\
&\quad+\langle \|u_0';H^{s+1}\| \rangle \|(u_0-u_0',\A_0-\A_0',\A_1-\A_1');X^{s-1,\sigma-\varepsilon}\|\},
\label{Dinitial}
\end{align}
where $D, E$ are replaced by $\tilD, \tilE$ if $s>2$.
Let $\eta$ be a rapidly decaying function such that $\eta$ is radial and $(\mathcal{F}\eta)(\xi)=1$ for $|\xi|\le 1$.
Here $\mathcal{F}$ is the Fourier transform.
We put $\eta_\delta(x)\equiv \delta^{-3}\eta(\delta^{-1}x)$ for $\delta>0$.
Then for $\theta$ with $0<\theta<\infty$, $\eta_\delta$ has the properties 
\begin{equation}
\|\eta_\delta\ast w;H^{s+\theta}\|\lesssim \delta^{-\theta}\|w;H^s\| \text{ for any } \delta>0,\
\|\eta_\delta\ast w-w;H^{s-\theta}\|= o(\delta^\theta) \text{ as } \delta\to 0.
\label{eta}
\end{equation}
Now we put 
\begin{align*}
u_0^\delta\equiv \eta_\delta\ast u_0, && 
\A_0^\delta\equiv \eta_{\delta^{1/\varepsilon}}\ast \A_0, && 
\A_1^\delta\equiv \eta_{\delta^{1/\varepsilon}}\ast \A_1,  \\
{u'_0}^\delta\equiv \eta_\delta\ast u_0', && 
{\A'_0}^\delta\equiv \eta_{\delta^{1/\varepsilon}}\ast \A_0', && 
{\A'_1}^\delta\equiv \eta_{\delta^{1/\varepsilon}}\ast \A_1',  
\end{align*}
and let $(u^\delta,\A^\delta)$ and $({u'}^\delta,{\A'}^\delta)$ be the solutions for the data 
$(u_0^\delta,\A_0^\delta,\A_1^\delta)$ and $({u'_0}^\delta,{\A'_0}^\delta,{\A'_1}^\delta)$, respectively. 
For fixed $\delta$, we have 
$D^{s,\sigma}(u^\delta,\A^\delta,{u'}^\delta,{\A'}^\delta)\to 0$ 
as $(u_0',\A_0',\A_1')\to (u_0,\A_0,\A_1)$ by virtue of \eqref{Dinitial}.
Therefore
\[
\limsup_{\rm data} D^{s,\sigma}(u,\A,u',\A')
\le D^{s,\sigma}(u,\A,u^\delta,\A^\delta)
+\limsup_{\rm data}D^{s,\sigma}(u',\A',{u'}^\delta,{\A'}^\delta).
\]
Here $\limsup_{\rm data}$ is the abbreviation of $\limsup_{(u_0',\A_0',\A_1')\to (u_0,\A_0,\A_1)}$.
The right-hand side is bounded by
\begin{align*}
&C(L)\{\|(u_0-u_0^\delta,\A_0-\A_0^\delta,\A_1-\A_1^\delta);X^{s,\sigma}\| \\ 
&\quad+\langle \|u_0^\delta;H^{s+1}\| \rangle \|(u_0-u_0^\delta,\A_0-\A_0^\delta,\A_1-\A_1^\delta);X^{s-1,\sigma-\varepsilon}\|\}.
\end{align*}
By \eqref{eta}, this tends to $0$ as $\delta \to 0$.
Therefore $\limsup_{\rm data} D^{s,\sigma}(u,\A,u',\A')=0$.
Repeating the process above finite times, 
we can show the result
on any compact interval where the unique existence of the solution is established.
\ebox

%%%%%%%%%%%%%%%%%%%%%%%%%%%%%% Lorentz and Temporal gauge
\section{The Lorentz gauge and the temporal gauge}

We prove only Theorem \ref{thmsl}; we can prove Theorem \ref{thmst} analogously.
In this section, we indicate by the superscript L the Lorentz gauge and by C the Coulomb gauge, respectively.
To begin with, we heuristically explain how we construct the solution to \textbf{MS-L}.
For any solution $(u^\LG,\phi^\LG,\A^\LG)$ to \textbf{MS-L}, there exists a solution to \textbf{MS-C}
which is gauge equivalent to $(u^\LG,\phi^\LG,\A^\LG)$.
Indeed, let us put $\lambda=\Delta^{-1}\Div\A^\LG$, $\A^\CG=P\A^\LG$, $\phi^\CG=(-\Delta)^{-1}\rho(u^\LG)$ and
$u^\CG=e^{-i\lambda}u^\LG$.
Then $\A^\LG=\A^\CG+\nabla\lambda$ by definition, and $\phi^\LG=\phi^\CG-\partial_t \lambda$ by the Lorentz gauge condition
and by the equation for $\phi^\LG$.
Therefore $(u^\LG,\phi^\LG,\A^\LG)$ and $(u^\CG,\phi^\CG,\A^\CG)$ are connected by the relation \eqref{gauge}.
Clearly $\A^\CG$ satisfies the Coulomb gauge condition, 
and $(u^\CG,\A^\CG)$ satisfies \textbf{MS-C} since \textbf{MS} is gauge invariant.
Moreover, $\lambda$ must satisfy the wave equation
\begin{equation}\label{eq:lambda}
(\partial_t^2 -\Delta)\lambda =\partial_t\phi^\CG
\end{equation}
with data $\lambda_j=\Delta^{-1}\Div\A_j$, $j=0,1$. 
The subscripts 0 and 1 indicate the initial datum for $\lambda$ itself and its time derivative, respectively.
Therefore we can solve the Cauchy problem for \textbf{MS-L} as follows. 
First we solve \textbf{MS-C} with data $(u_0^\CG,\A_0^\CG,\A_1^\CG)=(e^{-i\lambda_0}u^\LG,P\A_0^\LG,P\A_1^\LG)$.
Next we solve \eqref{eq:lambda}.
Then we construct the solution $(u^\LG,\phi^\LG,\A^\LG)$ to \textbf{MS-L} by the gauge transform.
\thmskip
\Proofof{Theorem {\rm \ref{thmsl}}}
We define $(u_0^\CG,\A_0^\CG,\A_1^\CG)$ as above.
Clearly $(\A_0^\CG,\A_1^\CG) \in H^\sigma \oplus H^{\sigma-1}$ by the boundedness of $P$ on $H^\sigma$.
Moreover if $\sigma\ge s-1$, $u^\CG \in H^s$ since $\lambda_0 \in \hdot^1\cap \hdot^{\sigma+1}$.
Therefore $(u_0^\CG,\A_0^\CG,\A_1^\CG)\in X^{s,\sigma}$.
By Theorem \ref{thmsc}, there exists an interval $I=[0,T]$ such that \textbf{MS-C}
with data $(u_0^\CG,\A_0^\CG,\A_1^\CG)$ has a unique solution $(u^\CG,\A^\CG)$ with 
$(u^\CG,\A^\CG,\partial_t\A^\CG)\in C(I;X^{s,\sigma})$.
The function $\lambda$ is obtained by the propagator $K(t)=\sin t\omega /\omega$ and its time derivative 
$\Dot{K}(t)=\cos t\omega$ as

%Then $\lambda$ satisfies , where we have used the expression 
\begin{align*}
\lambda&= \dot{K}(t) \lambda_0+K(t)\lambda_1
+\int_0^t K(t-\tau)\partial_\tau(-\Delta)^{-1}|u^\CG|^2(\tau)d\tau \\
&= \dot{K}(t) \lambda_0+K(t)(\lambda_1  +\Delta^{-1}|u_0|^2)
+\int_0^t \dot{K} (t-\tau)(-\Delta)^{-1}|u^\CG|^2(\tau)d\tau.
\end{align*}
Here we have used the integration by parts.
By this expression, we find that $\partial_t^j\lambda\in C^j(I;\hdot^1\cap \hdot^{\sigma+1-j})$ for $j=0,1$.
We define $(u^\LG,\phi^\LG,\A^\LG)$ as above.
Clearly this satisfies \textbf{MS-L}, the initial condition and the Lorentz gauge condition.
We can check that $(u^\LG,\A^\LG) \in C^j(I;H^{s-2j}\cap H^{\sigma-j})$ for $j=0,1$.
Moreover $\phi^\LG \in C^j(I;H^{\sigma-j})$ for $j=0,1$. Indeed $\phi^\LG$ satisfies
\begin{align}\label{eq:phiL}
\phi^\LG \equiv\phi^\CG-\partial_t\lambda
=\dot{K}(t) \phi_0^\LG+K(t)\phi_1^\LG
+\int_0^t K(t-\tau)|u^\LG|^2(\tau)d\tau
\end{align}
by virtue of the condition
$\phi_1^\LG +\Div\A_0^\LG=\Delta\phi_0^\LG+|u_0^\LG|^2 +\Div\A_1^\LG=0$.
The right-hand side of \eqref{eq:phiL} belongs to the desired space under the assumption for $(s,\sigma)$.
Therefore 
$(u^\LG,\phi^\LG,\partial_t\phi^\LG,\A^\LG,\partial_t \A^\LG)\in C([0,T];Y^{s,\sigma})$.
%Let $\phi'$ be the solution of 
%\[
%(\partial_t^2-\Delta)\phi'=|u'|^2,\ \phi'(0)=\phi_0',\ \partial_t\phi'(0)=\phi_1'.
%\]
%Then $\phi'$ satisfies $\phi'\in C([0,T];H^\sigma)$ and $\phi'=\phi-\partial_t\lambda$.
%Therefore $(u',\A',\phi')$ is a solution of \textbf{MS}, and it satisfies the Lorentz gauge condition.
%
%For the temporal gauge, let $\phi$ be the solution of
%\[
%\partial_t\lambda=\phi,\ \lambda(0)=\nabla\Delta^{-1}\A_0'.
%\]
%Then $\lambda$ satisfies $\omega\lambda \in C([0,T];H^\sigma)$.
%Therefore $u',\A',\phi'$ defined by \eqref{gauge} satisfies 
%$(u',\A')\in C([0,T];H^s\oplus H^\sigma)$ and $\phi'=0$.
%The gauge condition \eqref{temporal} follows from the conversation of the electric charge $\partial_t|u|^2=-\nabla\J(u,\A)$.
%
The uniqueness and the continuous dependence on the data of solutions follow from the well-posedness for \textbf{MS-C} and
\eqref{eq:lambda}. 
\ebox

\section*{Acknowledgments}
T. Wada would like to express his gratitude to professor J. Ginibre
for his kind hospitality, encouragement and valuable discussion at l'universit\'e de Paris XI.
The authors would like to thank professors Ginibre and G. Velo for their comments, 
by which the result is refined.

\bigskip

\bigskip

\noindent 
Makoto NAKAMURA 
\\*
\indent{\small
\begin{tabular}{l}
Graduate School of Information Sciences \\
Tohoku University \\ 
Sendai 980-8579, Japan \\
E-mail: m-nakamu@math.is.tohoku.ac.jp
\end{tabular}
}
\bigskip

\noindent Takeshi WADA
\\*
\indent{\small
\begin{tabular}{l}
Department of Mathematics \\
Osaka University \\
Osaka 560-0043, Japan \\
E-mail: wada@math.sci.osaka-u.ac.jp
\end{tabular}
}

\end{document}